\documentclass[10pt]{article}
\usepackage{amsmath,amssymb,amsthm,amscd}
\numberwithin{equation}{section}

\def\Bl{\biggl(}
\def\Br{\biggr)}
\def\p{\partial}
\def\b{\bar}

\def\e{\eta}

\def\l{\lambda}

\def\cL{{\cal L}}

\def\cD{{\cal D}}

\def\cH{{\cal H}}

\def\Tilde{\widetilde}

\def\cD{\mathcal D}

\def\cF{{\mathcal F}}

\def\cH{{\mathcal H}}

\def\cK{{\mathcal K}}
\def\cL{{\mathcal L}}

\def\cR{{\mathcal R}}
\def\cS{{\mathcal S}}

\def\cX{{\mathcal X}}

\def\bE{{\mathbf E}}

\def\bN{{\mathbf N}}

\def\CC{{\mathbb C}}
\def\NN{{\mathbb N}}

\newtheorem{prop}{Proposition}[section]
\newtheorem{theo}[prop]{Theorem}
\newtheorem{lem}[prop]{Lemma}
\newtheorem{cor}[prop]{Corollary}
\newtheorem{rem}[prop]{Remark}
\newtheorem{ex}[prop]{Example}
\newtheorem{defi}[prop]{Definition}
\newtheorem{conj}[prop]{Conjecture/Question}

\def\begeq{\begin{equation}}
\def\endeq{\end{equation}}

\newcommand {\beg}{\begin{eqnarray*}}
\newcommand {\ee}{\end{eqnarray*}}

\def\and{\quad{\rm and}\quad}

\def\and{\quad{\rm and}\quad}

\def\bN{{\mathbf N}}

\let\lra=\longrightarrow

\def\mapright\#1{\,\smash{\mathop{\lra}\limits^{\#1}}\,}

\begin {document}
\bibliographystyle{plain}
\title{Space of K\"ahler metrics III---On the lower bound of the Calabi energy and geodesic distance}
\author{Xiuxiong Chen\footnote{The author is partially supported by NSF grant. }\\ Department of Mathematics\\
University of Wisconsin,  Madison}
\date{}
% Enter your date or \today between curly braces
\maketitle
% 1. introduction
% 2. Known results ---cite results proved in previous paper but needed here.
% 3. Gap theorems
% 4. C^0 estimates
% 5. Re-adjust the automorphisms
% 6. C^2 and higher order derivatives estimates
% 7. The proof of main theorems.
\tableofcontents
\section{Introduction}
Inspired by the beautiful work of Donaldson \cite{Dona96}, the author initiated series of works \cite{chen991} \cite{chen992} aiming to understand  geometric structure of the space of K\"ahler potentials  and its application to interesting problems in K\"ahler geometry.  The existence of $C^{1,1}$ geodesic is established in \cite{chen991}, while the best regularity result on geodesic is given by
\cite{chentian005}.  However, the approach taken in \cite{chentian005}
is new (cf, Section 2 for further explanation.) and in many ways, the present work should be viewed as part III of this series.  It consists of three inter-related parts:
\begin{enumerate}
\item First, we prove a folklore conjecture on the greatest lower bound of the Calabi energy in any K\"ahler class.  This
was known in the 1990s when K\"ahler metrics is assumed to be invariant under some maximal compact subgroup \cite{Hwang951}.   
See ackwnoledgements for further remarks on this result.

\item Secondly, we give an upper/lower bound estimate of the K energy in terms of the geodesic distance and the Calabi energy.
This is used to prove a theorem on convergence of K\"ahler metrics in holomorphic coordinates, with uniform bound on the Ricci curvature and the diameter.  This kind of problems
is difficult because the K\"ahler geometry is more or less extrinsic while the well known Cheeger-Gromov
convergence theorem (with bound on curvature, diameters) is very
intrinsic.  
\item Thirdly, we set up a framework for the existence of geodesic rays when an asymptotic direction is given.  In particular,  if the initial geodesic ray is  tamed by a  bounded ambient geometry (c.f. Definition 3.2 and 3.3), 
then one  can derive some relative $C^{1,1}$ estimates for other geodesic rays in the same
direction.  More in depth discussions on  geodesic rays will be delayed to the beginning of Section 3.  In a sequel of this paper, we will give more regularity estimates on geodesic rays.
\end{enumerate}

\subsection{Donaldson's Conjectures}
According to Calabi \cite{calabi82},
an extremal K\"ahler metric is charicterized as the critical points of the $L^2$ norm of the
scalar curvature function in any given K\"ahler class.  The extremal K\"ahler metric includes the more
famous K\"ahler Einstein metric as a special case. 
In \cite{Dona96},  Donaldson set out an ambitious program in attacking core problems of  K\"ahler geometry via setting up formally a connection between geometric problems in the infinite dimesional
space and the current interesting problems in K\"ahler geometry.   In particular, he
  proposed  three inter-related conjectures:
\begin{enumerate}
\item The space of K\"ahler potentials is  uniquely connected by $C^\infty$ geodesic segments;
\item The space of K\"ahler potentials is a metric space;
% \item  The uniqueness of extremal K\"ahler metric follows from the smoothness of geodeisc;
\item  The non-existence of constant scalar curvature is equivalent to a geodesic ray where
the K energy functional(c.f. eq. (\ref{eq:kenergyform}) for definition) decays at $\infty.\;$
\end{enumerate}

   In \cite{chen991},  following Donaldson's program, the author established the existence of $C^{1,1}$
   geodesic by solving a Drichelet boundary value problem for a
   homogenous complex Monge Ampere equation.  Consequently, the second conjecture of Donaldson is completely verified.   Moreover, one important application is to show that Calabi's extremal
   K\"ahler metric (CextrK) is unique if the first Chern class is non-positive.    The uniqueness problem is completely settled now:  In algebraic manifold with discret automorphism group, 
   it is proved  by S. K. Donaldson \cite{Dona01}.   T. Mabuchi
 \cite{Ma05}  removes the assumption on the  automorphism group while X. X. Chen-G. Tian \cite{chentian005} complete the proof for general K\"ahler class.\\

   Chen-Tian \cite{chentian005} showed that the solution to disc version of geodesic
   problem is smooth except at most a codimension 2 set with respect to  generic boundary data.    For the convenience of readers, we will briefly describe the viewpoint of \cite{chentian005}
   in Section 2.   In particular, the partial regularity theory established in
\cite{chentian005} for solution of the disc version geodesic
equation plays a crucial role in this paper (for obtaining a lower bound of the Calabi energy).

   %As an application,  we  proved the uniqueness of extremal metric when the
   %first Chern class is non-negative.   These results already demonstrates that Donaldson's Program is indeed
   %very promising.    If regularity of geodesic can be improved to $C^{2,\alpha}$ or
   %even $C^2$ continuous,  then one can show that the K energy is weakly convex along geodesic
   %and the uniqueness of extremal K\"ahler metric in general  follows. \\
\subsection{Yau-Tian-Donaldson conjecture}
The Calabi conjecture on the existence of K\"ahler Einstein metrics has driven the subject for the second
half of the last century.   In late 1990s,
S. T. Yau conjectured that  the existence of K\"ahler
Einstein metric in Fano manifolds is equivalent to some form of
Stability of the underlying polarized K\"ahler class.  According to
G. Tian \cite{tian97} and Donaldson \cite{Dona96},  this equivalence relation
should be extended to  include the case of the constant scalar curvature (cscK) metric in a general
K\"ahler class. % Following \cite{}, we call this  ``Yau-Tian-Donaldson conjecture."  
   In a foundational paper,  G. Tian \cite{tian97} introduced the 
 notion of K Stability and in the same paper, he proved
 that the existence of KE metric implies weak K stability.    
 More recently,  in a fundamental paper \cite{Dona01},
Donaldson proved that, in algebraic manifold with discrete automorphism group, the existence of cscK
metric implies that the underlying K\"ahler class is Chow-Stable.  In this paper,  Donladson actually
formulated a new version (but equivalent) of K stability in terms of weights of Hilbert points. In K\"ahler 
 toric varieties, the existence of cscK metric implies the underlying K\"ahler class  is Semi-K stable \cite{Dona051}.    In  \cite{chentian005} , Chen-Tian proved that the existence
of a cscK metric implies that the K energy has a lower bound in this K\"ahler class.  Following the work of Paul-Tian \cite{paultian04}, this in turns implies the Semi-K stability of the underlying complex structure.  After we announced
our work \cite{chentian005},  S. K. Donaldson\cite{Dona05} proved a similar lower
bound in the algebraic settings.\\

\subsection{On the existence of geodesic rays}

The result of Donaldson on stability was extended by T.Mabuchi to the case of extremal
K\"ahler metric with some modified notion of stability. However,
for general K\"ahler classes, the usual notion of stability doesn't
apply because the manifold can not be embedded in $\CC P^N$ for
some large $N\gg 1.\;$ In \cite{Dona96}, Donaldson envisioned that
a geodesic segment or geodesic ray should play a similar role that
 a one parameter subgroup plays in a projective K\"ahler
manifold.   In the third conjecture of Donaldson's program,  he defines a set
of equivalence relationships:
\begin{enumerate}\item There exists no constant scalar curvature
metric in $(M,[\omega_0])$; \item There exists a geodesic ray from
some $\varphi_0\in \cH $ such that the K energy function
is strictly decreasing as $t \rightarrow \infty;\;$ \item From any
$\varphi \in \cH, $ there exists a geodesic ray initiated from
$\varphi$ such that the K energy function is strictly decreasing as
$t \rightarrow \infty.\;$
\end{enumerate}

The first step towards proving this conjecture is to establish an
existence result of geodesic ray with respect to another given
geodesic rays.  According to Calabi-Chen\cite{chen992}, the infinite
dimensional space $\cH$ is a non-positively curved space. By the 
Triangle comparison theorem, we can show that there always exists
a geodesic ray initiated from a given potential function in the direction of
any given geodesic ray.  However, the geodesic ray arisen this
fashion acquired very little regularity and it is very hard to use in practice.
As a first step in this direction, we prove 
\begin{theo} \label{geodesicray0} (cf. Theorem \ref{th:geodesicDegenerate})  If there exists a geodesic  ray(c.f.Defi. \ref{def:specialgeodesicray}) $\rho(t): [0,\infty) \rightarrow \cH$ which is tamed by a bounded ambient geometry, then for any K\"ahler potential $\varphi_0\in \cH$,  there exists a  relative $C^{1,1}$ geodesic ray $\varphi(t)$ initiated from
$\varphi_0.\;$  Moreover, this geodesic ray is parallel
to the original geodesic ray.
\end{theo}

A geodesic ray tamed by a bounded ambient geometry
is more or less ``parallel" to the notion of special degeneration of complex structure in
algebraic case.   In Section 3, we will discuss in length, various issues related 
to stability (in terms of geodesic rays).   It is expected that these notions are more or less equaivlent
to the corresponding notions in the algebraic settings.   We defer our discussions of this topic 
to the beginning of Section 3.

\subsection{On the lower bound of  geodesic distance and the collapsing of K\"ahler manifold}
The famous work of Cheeger-Gromov states that the set of Riemannian metrics with the following
three conditions:
\begin{enumerate}\item  uniform curvature bound, \item diameter is bounded from above,
\item volume is bounded from below,
\end{enumerate}
then this is compact under  $C^{1,\alpha} $ diffeomphism for some $\alpha \in (0,1).\;$ If the third condition is dropped,  then ``collpasing" may occur (volume converges to $0$).     On the other hand,
for any sequence  of K\"ahler metrics in a fixed K\"ahler class, the volume is {\it a priori}  fixed. 
With uniform control of the curvature and diameter from above,  it
will converge by subsequence to some K\"ahler metric with  perhaps a different complex structure.
In general,  we don't know what additional geometrical condition is needed to ensure that the limit complex
structure is the same with the original complex structure.     In fact, such  a sequence might collapse in
some Zariski open subset of the original K\"ahler manifold (i.e., the volume form vanishes
in this subset) while the subsequence of K\"ahler metrics converges as Riemmanian metrics up
to diffeomorphism (cf. \cite{ruan90}).
%Moreover,  the most collapsing direction gives rise to  a meromorphic vector field in this Zariski open subset.  In general, this
%collpasing vector field can not be extended to be a global holomorphic vector field (otherwise,
%up to holomorphic automorphisms, the sequence of K\"ahler metrics will converges in the holomorphic category.). 
 In the discussion below,  we will
refer this  phenomenon as ``K\"ahler collapsing. "    \\

One intriguing and  challenging question is: when this ``K\"ahler
collpasing"  occurs, does the geodesic distance (in the space of K\"ahler metrics) necessary diverge to $\infty?$ \footnote{ For  a sequence of metrics mentioned above which does not converge in the original
complex structure,  one is expected to prove, via implicite function theory, that the 
geodesic distance (to some fixed K\"ahler metrics) must diverges to $\infty.\;$ } This in turns leads to another question: how do we estimate the lower bound of the geodesic distance?  For instance, if the diameter of a sequence of K\"ahler metrics in a given K\"ahler class diverges to 
$\infty$,  does the geodesic distance of this sequence of K\"ahler metrics also diverge to $\infty?\;$\\

We first prove  a theorem which links  the K energy,  the Calabi
energy and the geodesic diameter together.   The author believes that this theorem is very interesting
in its own right.
\begin{theo} Let $\varphi_0, \varphi_1$ are two arbitrary K\"ahler potentials in the same K\"ahler class.  Then, the following inequality holds
\begin{equation}
\bE (\varphi_1) - d(\varphi_0, \varphi_1) \cdot \sqrt{Ca(\varphi_1)} 
\leq \bE(\varphi_0).
\end{equation}
Here $d(\varphi_0, \varphi_1)$ is the geodesic distance in the space of K\"ahler potentials.
\end{theo}

 In other words: if geodesic distance and Calabi energy is bounded, so is the upper bound of  the K energy.   This is quite surprising  since we don't know how to control the K energy, even after one assumes the uniform bound of the Riemannian curvature.    On the other hand,  fixing  $\varphi_1$ and let $\varphi_0$  change,  this formula gives a lower bound estimate of the K energy in terms of geodesic distance as well.  Clearly, this inequality
 is a natural generalization of the theorem \cite{chentian005} that the K energy has a lower bound
 if there is a cscK metric.   In fact, we conjecture that, in a fixed K\"ahler class, if the infimum
 of the Calabi energy reaches $0$, then the K energy must have a lower bound.\\
 
   An immediate corollary is:
 \begin{cor}  Let $\varphi$ be a K\"ahler potential such that its Calabi energy
 is bounded.   If $|\varphi|_\infty$ is bounded, then its K energy is bounded from above.
 \end{cor}

We say the K energy functional is ``proper"
if it is bounded from below by certain norm function which will be introduced in Section 2. 
We say the K energy functional is ``quasi-proper" if the K energy functional is bounded below by its highest order leading term (cf.  Section 2).  In K\"ahler Einstein manifold,  the K energy functional is always proper\cite{tian97}.   In a general K\"ahler manifold,  Tian conjectured that the cscK metric exists if and only if the K energy functional
is proper.    When the first Chern class is semi-negative,  there is a sufficient condition that  the K energy functional in that K\"ahler class is either proper or quasi proper \footnote{For instance, in complex dimension 2, if
\[
 {{ [\omega] \cdot [ - C_1(M)]}\over { [\omega]^{[2]}}} (- C_1(M)) - [\omega] > 0 
\] 
then the K energy is quasi-proper \cite{chen00}.   For higher dimension K\"ahler manifolds, readers are referred to
Song-WeinKove\cite{SongBen040} .}.\\

Now we are ready to answer the question about ``K\"ahler collpasing": 
 \begin{theo}  (No ``K\"ahler collapsing") Let $(M, [\omega])$ be a polarized K\"ahler manifold where the K energy functional
 is either proper or quasi-proper.    Let $\cS$ be a set of K\"ahler metrics in $[\omega]$ 
with uniform Ricci curvature bound from below and diameter bound from above \footnote{The diameter bound can be replaced by a bound on Sobleve constant.}.  If this set of K\"ahler metrics lies in a bounded geodesic ball
in the space of K\"ahler metrics, and  if we assume Ricci also has an upper bound, then
all metrics in $\cS$ are uniformly equivalent to each other in $C^{1,\alpha}(M)$ topology for any $\alpha \in (0,1).\;$
In particular, ``K\"ahler collapsing" will not occur.
\end{theo}

Note that the geodesic distance appears to be a very weak notion.   The bound on Ricci curvature is much weaker than the conditions stated in Cheeger-Gromov's theorem.  However, 
the combination of the two conditions seems to be very powerful.   \\

In a subsequent paper,  we will drop the assumption that the Ricci is bounded from above. The assumption that the K energy functional is either proper or quasi-proper in $(M, [\omega])$ is just technical.   We hope that this will be removed in a subsequent work.  
% \begin{theo}  \footnote{I DON'T have  a proof of this theorem yet and this is not a final version either.  It may  have to wait for next paper. }Let $(M, [\omega])$ be a polarized K\"ahler manifold where either  the K energy functional  is proper or the first Chern class is negative.    Let $\cS$ be a set of K\"ahler metrics in $[\omega]$  with Ricci curvature bounded from below, uniform bound on Calabi energy and Sobleve constant.  If this set of K\"ahler metrics lie in a bounded geodesic ball in the space of K\"ahler metrics, then all metrics (by subsequence) in $\cS$ are uniformly equivalent to each other in any compact subset of a open dense subset of $M, $ where the complementary of this dense open subset is a subvariety of Hausdorff codimension 2. \end{theo}

\subsection{On the lower bound of the Calabi energy } 
It is well known that the Calabi energy is locally convex near an extremal K\"ahler metric.  It is a very interesting and difficult question if the Calabi energy in the K\"ahler class is bounded below by the energy of the extremal metric. According to Calabi\cite{calabi82}, an extremal K\"ahler metric automatically exhibits the maximal symmetry possible allowed by the underlying complex structure.  In the 1990s,  A. Hwang \cite{Hwang951} proved that the Calabi energy of the invariant K\"ahler metrics (maximal possible symmetric...) is bounded below by the absolute value of the Futaki invariant (evaluated at the Canonical extremal vector field). If there is an extremal K\"ahler metric in this
class,  the absolute value of the Futaki invariant is precisely  the Calabi energy of the extremal
K\"ahler metrics.  Hwang's proof uses strongly the bi-invariant metric
in the Lie algebra of gradient holomorphic vector fields where the symmetric property is the key
to define this positive definite invariant metric.   In 1980s, it is conjectured that the same low bound
holds for all metrics in the same K\"ahler class.   There are many attempts to generalize this to all
K\"ahler metrics,  and this problem has proved to be very difficult indeed.\\

Aside from this Folklore conjecture on the Calabi energy, there are other important motivations to
study the lower bound of the Calabi energy {\it a priori}, for instance, the issue related to stability
and degeneration of K\"ahler manifolds.  For our strategy to work, the  main technical obstacle has been the insufficient regularity of
the $C^{1,1}$ geodesic. 
However, the partial regularity theory established in
\cite{chentian005} for solution of the disc version geodesic
equation plays crucial role here. In particular, it is a very
powerful fact that the restriction of the K energy of this family
of K\"ahler potential over disc is subharmonic. We are able to use
this fact to
establish a lower bound for the Calabi energy in terms of
any effective destabilizing geodesic ray (cf.
Defi.\ref{def:geodesicdestable}).  In  particular, we prove this folklore conjecture about the
lower bound of the Calabi energy in each K\"ahler class. \\

\begin{theo}\label{th:calabienergybound3} Let $\cK$ be the Lie algebra of complex gradient holomophic vector field.   Then for any K\"ahler metric $\omega_g$ in $[\omega]$, we have

\[
  Ca(\omega_g) \geq  \cF_{\cX_c}([\omega]) .
\]
where $\cX_c$ is the {\it a priori} extremal vector field in $(M, [\omega])\;$ and $\cF$ is the Futaki invariant. The equality holds when
$g$ is an extremal K\"ahler metric.
\end{theo}
More generally, we have

\begin{theo}\label{th:calabienergybound0} Suppose $\rho(t):[0,\infty) \rightarrow \cH$ is an effective
destabilized geodesic ray
in $\cH$,  then
\begin{equation} \displaystyle \inf_{\varphi
\in\cH} \int_M\; (R(\varphi)- \underline{R})^2 \omega_\varphi^n
\geq \displaystyle \sup_{\rho} \yen(\rho)^2,
\label{th:calabienergybound00}
\end{equation}
where the sup in the right hand side 
runs over all possible destabilized geodesic rays.   The definition of an effective destabilized geodesic ray and $\yen$ invariant are given in Defi.\ref{def:geodesicdestable} and Defi. \ref{defi:geodesicinvariant} respectively.
\end{theo}
\begin{defi}
Let $(M, [\omega_0], J_0)$ be a triple K\"ahler structure.  Another triple K\"ahler structure
$(M',[\omega'], J')$ lies in the closure of the diffeomorphism orbit of $(M, [\omega_0], J_0)$ if there exists a sequence of K\"ahler forms $\{\omega_{\varphi_m}, m\in \NN\} \subset [\omega]$  and a sequence of diffeomorphism 
$\{f_m \in Diff(M), m \in \NN\}$ such that $ (M, f_m^*\omega_{\varphi_m}, f_m^*J_0)$ converges to $(M', \omega', J')$ in $C^{1,\alpha}$ topology for some $\alpha \in (0,1).\;$
\end{defi}

\begin{defi}  Let $(M, [\omega_0], J_0)$ be a triple K\"ahler structure.  Suppose that $(M',[\omega'], J')$
is another   triple K\"ahler structure
 which lies in the closure of the diffeomorphism orbit of $(M, [\omega_0], J_0).\;$ $(M',[\omega'], J')$ is called destabilizer of the original triple K\"ahler structure $(M, [\omega_0], J_0)$ if  there exists an effective destabilized geodesic ray $\varphi(t)$   in $(M, [\omega_0], J_0)$ such 
 that   there is a subsequence
 of $\omega_{\varphi(t)} (t\rightarrow \infty)$ which converges to a metric in $(M', [\omega'], J')$
 up to diffiomorephism.
\end{defi}

Now we can extend Theorem \ref{th:calabienergybound3} 
to a more general setting:
\begin{theo}\label{th:calabienergybound4} Let $(M, [\omega_0], J_0)$ be a triple K\"ahler structure.
%  Suppose $\b \cF$ is ths supremum of $\cF_{\cX_c}([\omega'])$ over all K\"ahler triple structures
%$(M', [\omega'], J')$  in the closure of the diffeomorphism orbit of $(M, [\omega_0], J_0).\;$  
The following inequality hold
 \[
 \displaystyle \inf_{g \in \cH}  Ca(\omega_g) \geq \displaystyle \sup_{(M',[\omega'], J')}  \displaystyle \inf_{ (X, X)=1, X \in \cK(J')}  {(X, \cX_c)^2\over (X,X)}.
\]
Here the inner product is the Futaki-Mabuchi inner product for the Lie algebra $\cK(J')$ of the Maximum compact subgroup $K(J')$ of $Aut(M', J').\;$   The supremum runs over all possible K\"ahler triple structures
$(M', [\omega'], J')$  in the closure of the diffeomorphism orbit of $(M, [\omega_0], J_0)\;$ which
destabilized  $(M, [\omega_0], J_0)$. \end{theo}

One should be able to define a weak notion of destablizing K\"ahler triple later, while the inequality
in the preceding theorem still hold.

\begin{defi}  Suppose the K\"ahler triple $(M, \omega, J)$ satisfies the following inequality
\begin{equation}
\cF_{\cX_c}([\omega_0]) \geq \displaystyle \sup_{(M',[\omega'], J')}  \cF_{\cX_c}([\omega']),
\label{eq:diffgeomstable}
\end{equation}
where $(M', [\omega'], J')$ are all K\"ahler triples in its closure of diffemorphisms.   Then we
call $(M, [\omega], J)$ stable in the sense of differential geometry.
\end{defi}
An immediate  interesting/challenging question is: what is the relation of stablity in the sense
of differential geometry with other notions of stability such as K stablity? In algebraic manifold,
the notion of K stability shall be stronger than this one.\\

In light of these theorems, one expects that there is a deep
relation between geodesic rays, test configurations and their
respective role in defining stability.  We then propose some
notions of stability in terms of geodesic rays which might be viewed
as a natural extension of what is given in \cite{Dona96}.  Moreover, the
relation between geodesic stability and K stability should be also
an interesting topic to explore in near future.  More extensive discussions
on this topic will be delayed to Section 3.\\
  
\noindent {\bf Organization:} In Section 2, we
give a brief outline of known results in the space of K\"ahler
potentials.    In
Section 3, we prove that the existence of geodesic ray with respect to some 
``nice" geodesic ray.  In
Section 4, we give a greatest lower bound estimate for the Calabi energy.
In Section 5, we give a lower bound estimate of the geodesic distance and
rule out the possibility of ``K\"ahler collapsing" in bounded geodesic balls in $\cH.\;$ \\
%In Section 6,  we propose some conjectures in K\"ahler geometry. \\

\noindent{\bf Acknowledgment}:   The strategy of obtaining a lower
bound of the Calabi energy through geodesic ray has been discussed
with S. K. Donaldson in 1997-98,   on and off  since then.  The author wants to thank Professor
Donaldson for kindly sharing his insight on this matter.  Readers are encouraged to
compare the results on lower bound of Calabi energy to \cite{Dona051} (In particular,
Theorem 1.6.). \\

      Thanks also goes to Professor G. Tian,  our joint work \cite{chentian005} really provided
      technical support to Theorem 1.2.    Thanks also to Professor Calabi  and Professor J. P. Bourguinon for their continuous support in my research in last few years.  My student
Yudong Tang  carefully read through an earlier version of this paper and I want to thank him
for his help.      

\section{Brief outline of geometry in the Space of K\"ahler potentials.}
\subsection{Quick introduction of K\"ahler geometry}
Let $\omega$ be a fixed K\"ahler metric
on $M$. In a holomorphic coordinate, $\omega$ can be expressed as
\[
\omega = g_{\alpha\b \beta} \;{\sqrt{-1}\over 2}\;d w^\alpha \wedge
d\,w^{\b \beta}.
\]
The Ricci curvature can be conveniently expressed as \[
R_{\alpha\b\beta} = -  {{\p^2 \log \det \left(g_{i\b j}\right)
}\over {\p w^\alpha \p w^{\b \beta}}}.
\]
The scalar curvature can be defined as
\[
R  = - g^{\alpha\b \beta }\;  {{\p^2 \log \det \left(g_{i\b
j}\right) }\over {\p w^\alpha \p w^{\b \beta}}}.
\]
The so called Calabi energy is
\begin{equation}
Ca(\omega) = \int_M \; (R(\omega) - \b R)^2 \omega^n.\;
\label{eq:calabienergy}
\end{equation}
Here $\b R$ is the average scalar curvature value for all metric
in the K\"ahler class.  According to Calabi 
\cite{calabi82}\;\cite{calabi85}, a K\"ahler metric is called extremal if the
complex gradient vector field
\begin{equation}
\cX_c = g^{\alpha\b \beta} {{\p R}\over {\p w^{\b \beta}}}
{\p\over {\p w^\alpha}} \label{eq:extremalvectorfield}
\end{equation}
is a holomorphic vector field.  According to \cite{FutakiMa95}, the
extremal vector field $\cX_c$ is {\it a priori }  determined in
each K\"ahler class, up to holomorphic conjugation.\\

If $X$ is a holomorphic vector field, then for any K\"ahler potential
$\varphi$ we can define $\theta_X$ up to some additive constants by
\begin{equation}
L_X \omega_\varphi = \sqrt{-1} \p \b \p \theta_X.
\label{eq:liederivatives}
\end{equation}
Then, the well known Calabi-Futaki invariant
\cite{Futaki83}\;\cite{calabi85} is 
\begin{equation}
\cF_X([\omega]) =  \displaystyle
\int_M\; \theta_X \cdot (\underline{R} - R(\varphi))
\;\omega_\varphi^n. \label{eq:calabifutakiinvariant}
\end{equation} Note that this is a Lie algebra
character which depends on the K\"ahler class only.

\subsection{Weil-Peterson type metric by Mabuchi}
It follows from the Hodge
theory that the space of K\"ahler metrics with K\"ahler class
$[\omega]$ can be identified with the space of K\"ahler potentials
\[
{\cH}= \{ \varphi \mid \omega_{\varphi} = \omega + \b \p
\p \varphi > 0,\;{\rm on} \; M\}/ \sim,
\]
where $\varphi_1 \sim \varphi_2$ if and only if
$\varphi_1=\varphi_2+ c$ for some constant $c$. A tangent
vector in $T_\varphi\cH$ is just a function $\psi$ such that
$$\int_M \psi\omega_\varphi^n =0.$$ Its norm in the $L^2$-metric
on $\cH$ is given by (cf. \cite{Ma87})
\[
\|\psi\|^2_{\varphi} =\int_{M}\psi^2\;\omega^n_\varphi.
\]
It was subsequently defined similarly in \cite{Semmes92} and
\cite{Dona96}.    In all three papers, \cite{Ma87}\cite{Semmes92} and \cite{Dona96},  the authors
defined  this Weil-Peterson type metric from various points of view and proved formally that this infinite
dimensional space has non-positive curvature.   Using this definition,  we can define a distance
function in $\cH$:  For any two K\"ahler potentials $\varphi_0,
\varphi_1\in \cH$, let $d(\varphi_0,\varphi_1)$ be the infimum of
the length of all possible curves in $\cH$ which connects $\varphi_0$
with $\varphi_1.\;$  \\

A straightforward computation shows that a geodesic path
$\varphi:[0,1]\rightarrow \cH$ of this $L^2$ metric must satisfies
the following equation
\[
\varphi''(t) - {g_\varphi}^{\alpha \b \beta} {{\p^2 \varphi}\over
{\p t \p w^\alpha}}  {{\p^2 \varphi}\over {\p t \p w^{\b \beta}}}
= 0.
\]
where
\[
g_{\varphi,\alpha\b \beta} =  g_{\alpha \b \beta} + {{\p^2
\varphi}\over {\p w^\alpha\p w^{\b \beta}}}.
\]
 
According to S. Semmes \cite{Semmes92},   this path $\{\varphi(t)\}$ satisfies the
geodesic equation if and only if the function $\phi$ on $[0,1]\times
S^1\times M$ satisfies the homogeneous complex Monge-Ampere equation
\begin{equation}
\label{eq:hcma0} (\pi_2^* \omega + \partial \overline{\partial
}\phi)^{n+1} \;\; =\;\; 0, \qquad {\rm on} \; \Sigma \times M,
\end{equation}
where $\Sigma = [0,1] \times S^1$ and $\pi_2: \Sigma\times M\mapsto
M$ is the projection.  In fact, one can consider (\ref{eq:hcma0})
over a general Riemann surface $\Sigma$ with boundary condition
$\phi =\phi_0$ along $\partial \Sigma$, where $\phi_0$ is a smooth
function on $\partial \Sigma \times M$ such that $\phi_0(z,\cdot)\in
\cH$ for each $z\in \partial \Sigma$.\footnote{We often regard
$\phi_0$ as a smooth map from $\p \Sigma$ into $\cH$.} It also has
geometric meaning. The equation (\ref{eq:hcma0}) can be regarded as
the infinite dimensional version of the WZW equation for maps from
$\Sigma$ into $\cH$ (cf. \cite{Dona96}).\footnote{The original WZW
equation is for maps from a Riemann surface into a Lie group.}
\\

Next we introduce three  well known  functionals in $\cH$
here.   First, the so called $I$ functional is defined as
\[
{{d I(\varphi(t))}\over {d\, t}} =  \int_M\; {{\p\,\varphi}\over {\p\,t}} \omega_{\varphi(t)}^n, \qquad \varphi(t) \in \cH.
\]
The advantage of the and $I$ functional is that it is a constant along geodesic.   One can write down
an explicit formula for $I$ functional
\begin{equation}
I(\varphi) = \int_M\; \varphi \;\omega^n - \sum_{k=0}^{n-1}\; {{n-k}\over {n+1}}  \;\p \varphi \wedge \b \p \varphi\wedge \omega^k \wedge \omega_\varphi^{n-k-1}.
\label{eq:Ifunctional1}
\end{equation}
We write down a detailed proof for $I(\varphi)$ in Section 5
(prop \ref{prop:Ifunctional}).\\

Secondly, the so called $J$ funtional is defined as
\[
J(\varphi) = \int_M\; \varphi \left(\omega^n -\omega_\varphi^n\right) = 
\int_M\; \p \varphi \wedge \b \p \varphi \wedge \left(\displaystyle \sum_{k=0}^{n-1} \omega^k \wedge \omega_\varphi^{n-k-1} \right) > 0.
\]
Finally, the 
 K energy functional
(introduced by T. Mabuchi) is defined as a closed form $d\,\bE$.
Namely, for any $\psi \in T_\varphi \cH$, we have
\begin{equation}
(d\,\bE, \psi)_\varphi = \int_M\; \psi\cdot (\underline{R} -
R(\varphi)) \;\omega_\varphi^n.\label{eq:kenergyform}
\end{equation}
Note that for any holomorphic vector field, we have
\[
\cF_X([\omega]) = (d\,\bE, \theta_X)_\varphi.
\]
The K energy functional is called proper in $(M, [\omega]) $ if there exists a small constant $\delta > 0$
a constant $C$ such that
\[
  \bE (\varphi) \geq J(\varphi)^\delta - C.
\]
The K energy functional is called proper in $(M, [\omega]) $ if there exists a small constant $\delta > 0$
and a  constant $C$ such that
\begin{equation}
  \bE (\varphi) \geq \delta \int_M \log {\omega_\varphi^n\over \omega^n} \omega_\varphi^n  -  C
  \label{eq:quasiproper}
\end{equation}

\subsection{The new approach in Chen-Tian's paper
\cite{chentian005}}
The lack of sufficient regularity in Chen's solution to geodesic equation \cite{chen991} is the
obstruction to proving deeper  results in general cases by
using geodesics approaches.   We should point out that complex
Monge-Ampere equations have been studied extensively (cf.
\cite{CNS84}, \cite{CKNS85}, \cite{Bedford76} etc.). However,
regularity for solutions of homogeneous complex Monge-Ampere
equations beyond $C^{1,1}$ has been missing. Indeed, there are
examples in which some solutions are only $C^{1,1}.$\\

    The best regularity result about geodesic segment is
   due to Chen-Tian \cite{chentian005} where they showed that the solution to the disc version of the geodesic
   problem is smooth except at most a codimension 2 set with respect to  generic boundary data.   Although the $C^{1,1}$ bound
   derived in \cite{chen991} plays a crucial role,   Chen-Tian takes a new viewpoint towards geodesic
   equation.     The main theorem in \cite{chentian005} is 
\begin{theo}\cite{chentian005}
\label{th:almostsmooth} Suppose that $\Sigma$ is a unit disc. For
any $C^{k,\alpha}$ map $\phi_0: \p \Sigma \rightarrow \cH$ ($k \ge
2$, $0< \alpha < 1$) and for any $\epsilon
> 0$, there exists a $\phi_\epsilon: \p\Sigma\rightarrow \cH$
in the $\epsilon$-neighborhood of $\phi_0$ in
$C^{k,\alpha}(\Sigma\times M)$-norm, such that (\ref{eq:hcma0})
has an almost smooth solution with boundary value $\phi_\epsilon$.
\end{theo}

An almost smooth solution of eq. \ref{eq:hcma0} has a uniform
$C^{1,1}$ bound and smooth almost everywhere.  A detailed explanation
(including definitions) can be found in \cite{chentian005}.  However,
the importance of this theorm lies in the following

\begin{theo} \cite{chentian005}
\label{th:hessianofK-energy} Suppose that $\phi$ is a partially
smooth solution to (\ref{eq:hcma0}). For every point $z \in
\Sigma, $, let $\bE(z)$ be the K-energy (or modified K energy)
evaluated on $\phi(z,\cdot)\in \overline{\cH}$. Then $\bE$ is a
bounded subharmonic function on $\Sigma$ in the sense of
distributions, moreover, we have the following
\[\int_{\cR_\phi} | {\cD}
{{\p \phi}\over {\p z}}|_{{\omega_{\phi(z,\cdot)}}}^2\,
{\omega_{\phi(z,\cdot)}}^n\, d\,z  d\, \b z \leq \displaystyle
\int_{\p \Sigma}\;{{\p \bE}\over{\p {\bf n}}} \big |_{\p \Sigma}
ds,
\]
where $ds$ is the length element of $\p \Sigma$ and for any smooth
function $\theta$, $\cD \theta$ denotes the (2,0)-part of
$\theta$'s Hessian with respect to the metric
$\omega_{\phi(z,\cdot)}$. The equality holds if $\phi$ is almost
smooth.
\end{theo}

\section{On the existence of geodesic ray}
\subsection{Definitions and main results} 
  
  As suggested in \cite{Dona96}, one should view the geodesic rays as effective substitute
of a one parameter family of
subgroup acting on projective K\"ahler manifolds. % (c.f. \cite{Phong_Sturm05} also).  
   It is natural to compare geodesic rays to the test configurations $\pi: \cX \rightarrow \triangle$ such that all fibre K\"ahler manifolds $\pi^{-1}(t)$ are 
bi-holomorphic to each other  except when $t= 0.\;$  The central fibre usually carries a different complex structure with singularities.  The generic case is so called ``Normal cross singularity" and the special
case is when the central fibre is either smooth or has singular local codimenison 4 or higher.
However, by blowing up a few points in the central fibre if necessary, it might be possible to make the
total space smooth or  have  some form of bounded geometry.   For any test configuration, it might be possible to prove that 
there is always a relatively $C^{1,1}$ geodesic ray which is asymtotically closed to the test configuration
near the central fiber.  If the central fibre is smooth or smooth except a subvariety of condimension 4,
then the geodesic ray is smooth  generically except perhaps a singular locus of codimenion
two or higher.  \\

Motivated from the study of test configuration in algebraic setting,    in this section, we
 restrict our attentions to the case of ``nice"  geodesic rays $\omega_{\rho(t)} (t\in [0,\infty))$ which satisifies the following conditions:
 \begin{enumerate}
 \item  The non-compact family $ \left( M,  \omega_{\rho(t)} \right)$
 can be  compactified in some sense;
%\item   The curvature of $ \pi_2^* \omega_0 + \omega_{\rho(t)}  $ is uniformly bounded
 %from below;
 \item The limit of $(M, \omega_{\rho(t)})$ as $t\rightarrow \infty$ under suitable topology
 is smooth in the ``compactiftication" or has mild singularities (codimension 4 and higher). 
% \item   The curvature  of $(M, \omega_{\rho(t)})$ is uniformly bounded and the injectivity is uniformly
 %from below.
 \end{enumerate}

 The most special case of geodesic rays are those arising  from a fixed gradient complex holomorphic vector field.  In this case, the curvature  of $(M, \omega_{\rho(t)})$ is uniformly bounded and the injectivity radius is uniformly bounded from below.\\

Consider $\pi_2: ([0,\infty)\times S^1 ) \times M \rightarrow M$ as natural projection map.
\begin{defi}  A path $\rho(t):[0,\infty) \rightarrow \cH$ is called strictly convex if
$\pi_2^* \omega_0 + i \p \b \p \rho$ defines a K\"ahler metric in $([0,\infty)\times S^1 ) \times M.\;$
\end{defi} 

\begin{defi}\label{def:specialgeodesicray} A geodesic ray
$\rho(t) (t \in [0,\infty))$ is called special if it is one of the following types:

\begin{enumerate}
\item {\bf effective} if  the Calabi energy of $\omega_{\rho(t)}$ in $M$ is dominated by ${\epsilon\over t^2}$ for any $\epsilon > 0$ as $t \rightarrow \infty.\;$
%uniformly bounded for $t\in [0,\infty).\;$
\item{\bf  normal}  if the curvature of $\omega_{\rho(t)} $ in $M$ is uniformly bounded  for $t\in [0,\infty).\;$
\item {\bf bounded geometry} if  $(M, \omega_{\rho(t)}) (t\in [0,\infty)) $  has uniform bound on  curvature and a uniform positive lower bound of injective radius. 
\end{enumerate}
\end{defi}
\begin{defi} {\bf Tamed by a bounded ambient geometry}  A K\"ahler metric $h = \pi_2^* \omega_0 + i \p \b \p \b \rho$ 
in  $([0,\infty)\times S^1) \times M$ is said to have bounded ambient geometry if
\begin{enumerate}
\item  it has a uniform bound on its curvature;
\item  $([0,T]\times S^1 \times M, h)$ has a uniform lower bound on injectivity radius
 and the bound is independent of $T\rightarrow \infty;\;$
\item The vector length $|{\p\over {\p t}}|_h$ has a uniform upper bound. 
 \end{enumerate}
    A geodesic ray $(M, \omega_{\rho(t)})$
 is called tamed by this ambient metric $h$ if  there is a uniform bound of the relative potential
 $\rho-\b \rho$, or if there is a uniform constant $C$  such that 
 \footnote{It is possible to only assume these two inequsalities holds for a sequence of $t_i \rightarrow \infty$. }%The inequality
% implicitely assume that $|{\p \over {\p t}}|_h$ has a uniform lower bound.  Otherwise,
% we need to assume the two inequalities holds for any $(t,\cdot)$ such that the distance in $h$
 %from $(t,\cdot)$ to $\{t_i\}\times S^1\times M$ is finite (say less than $1.\;$). }
 \begin{enumerate}
\item $\displaystyle \max_{t} \;| n+1 +  \triangle_h (\rho - \b \rho)|  \leq C;$
\item $\displaystyle \max_{t}\;  |{{\p (\rho -\b \rho)}\over {\p t}} |_h \leq C;$
%\item $|{\p \over {\p t}}|_h \leq C.$
\end{enumerate}
%Here $h = \pi_2^* \omega_0 + i \p \b \p \b \rho.$
 %\end{enumerate}
\end{defi}

\begin{rem} A geodesic ray, tamed by a bounded ambient geometry,  corresponds to the special degeneration
of the complex structure in the algebraic setting. % From the constructions in the example below, we expect
%that these geodesic rays must also have bounded ambient geometry.  A geodesic
%ray with bounded ambient geometry corresponds to the case where the total moduli space is smooth in the test configuration.
 In the future, we should broad our definition of bounded ambient geometry to  include the following situations:
\begin{enumerate}
\item The upper bound of the curvature of the ambient K\"ahler metric might not be uniform;
\item The injectivity radius may have a lower bound which depends on the distance to some singular subvariety of higher codimension as well as on $t;\;$
\item   The restriction of the ambient K\"ahler metric $h$ in  $\{t_i\} \times S^1 \times M$ may have
some finite geodesic distance to $(M, \omega_{\rho (t_i)})\;$ while the later has certain geometric
bounds (such as the Calabi energy or  Sobelev constant, cf. Theorem 1.4).
\end{enumerate}
Of course, the regularity of geodesic may be weaked a bit as well. 

\end{rem}

\begin{rem} Using Cauchy-Kowalevaski's classical theorem,   Arezzo-Tian \cite{Are-Tian03}
proved that, a special degeneration of a complex structure when the
central fibre is analytic,  is asymptotically equivalent to
a geodesic ray  near the central fibre. 
\end{rem}
\begin{ex}  Suppose that $X$ is a gradient holomorphic vector field and let $\omega_0$ be a K\"ahler form invariant under $Im(X).\;$  Let $\sigma(t), t\in[0,\infty) $ be the automorphism group generated by $X.\;$
Set
\[
 \omega_{\rho(t)} = \sigma_t^* \omega_0.\]
 A straightforward calculation shows that $\rho(t) (t \in (-\infty, \infty)$ is a geodesic line.   Let $\sigma = \sigma_1$ and let $g_1 = \sigma^* g_0$ and $g_0$ be the two  K\"ahler metrics corresponding to $\omega_0$ and $\sigma^*\omega_0.\;$ Note that \[
 z {\p \over {\p z}} + X\]
  induces a $\CC^*$ action
$\b \sigma$ on $\Delta \times M\;$ which coincides with $\sigma $ in the manifold direction and
the multiplicity action on $\Delta$  direction.  Let $z_0=1$ and $z_k = \b \sigma^k z_0 \rightarrow 0.\;$ 
Set \[
M_{l,k} = \{ {1\over 2^l} \leq |z| \leq {1\over 2^k}\} \times M, \qquad \forall\;\; l,k \in \NN. 
\]
Then
\[
M_{0, \infty} = (\Delta\setminus \{0\}) \times M.
\]
It is easy to see that there is a smooth $S^1$ invariant  K\"ahler metric $\b h$ in $M_{0,1}$ such that 
\begin{enumerate} \item $\b h \mid_{|z|=0} = g_0$ and $\b h\mid_{|z|=1} = g_1.\;$
\item $h$ and $\b \sigma^* \b h$ gives rise a smooth metric in $M_{0,2}.\;$
\end{enumerate}  
Using $\b h$, we can define a K\"ahler $h$ in $ (\Delta\setminus \{0\}) \times M$  simply by
\[
h (z, \cdot)  = {\b \sigma}^{k*} \b h, \qquad \forall\;  (z, \cdot) \in M_{k, k+1}.\]
   By definition, $h$ is a smooth metric in $(\Delta\setminus \{0\}) \times M$ which has bounded curvature
and uniform positive lower bound on injectivity radius. \end{ex}

In fact, any normal geodesic ray is expected to be tamed by some bounded ambient
geometry, at least when it has bounded geometry.   

%\begin{rem} From the construction of the background metric with bounded geometry in the preceding
%example, all we need is a discrete $\CC^*$ action in $M$ which reflects the degeneracy of the given
%geodesic ray.  Since the limit of a geodesic ray with bounded geometry is a geodesic line arising from
%a gradient holomorphic vector field (perhaps in a different complex structure).  As in the preceding
%example, one can construct a background K\"ahler metric (in the new complex structure)
%with bounded geometry but reflecting the degeneracy of this geodesic line.  By perterbation argument,
%one should be able to construct a K\"ahler metric in the original complex structure with bounded geometry
%which reflects the degeneracy of the initial geodesic ray with bounded geometry.
%\end{rem}
\begin{defi}\label{def:geodesicparallel} For any two geodesic rays $\rho_1(t), \rho_2(t):
[0,\infty)\rightarrow \cH,\;$  they are called parallel if there
exists two constants $ C$ such that
\[
  \displaystyle \sup_{t \in [0,\infty)}\; \left(\rho_1(t)-\rho_2(t)\right)\leq C.
\]
\end{defi}

\begin{theo} \label{th:geodesicDegenerate} If there exists a  geodesic ray  $\rho(t): [0,\infty) \rightarrow \cH$ which is tamed by an ambient geometry,
then for any K\"ahler potential $\varphi_0\in \cH$,  there exists a relative $C^{1,1}$ geodesic ray $\varphi(t)$ initiated from
$\varphi_0$ and parallel to $\rho(t)$  such that \[
\displaystyle \sup_{t \in [0,\infty)}\; \mid\rho(t)-\varphi(t)
  \mid\leq C.
\] 
 Denote $\rho -\b \rho$ is the relative K\"ahler potential of the given geodesic ray
 with respect to the ambient metric $h.\;$  If $\rho-\b \rho$ is uniformly $C^2$ bounded
 in manifold direction and if it is $C^1$ bound in $t$ direction with respect
 to the ambient metric $h$,  then there exists two uniform constants $, \lambda,  C$ such that

\[
0\leq n+1 + \tilde\triangle (\varphi(t,x)-\rho(t,x))\leq C \exp {\lambda (\rho(t,x) - \b \rho(t,x))} .
\]
Here $\tilde\triangle$ is taken with respect to the ambient  K\"ahler
metric $h.\;$ In particular, when\footnote{This is the case when the geodesic ray
is given by one parameter holomorphic transfermation (c.f. ex. 3.5).} 
\[
\rho(t, x) - \b \rho(t, x)
\] is uniformly bounded, the resulting geodesic ray has a uniformly $C^{1,1}$  bound
in terms of ambient metric $h.\;$
The constant $\lambda, C$ depends on $h.\;$ \end{theo}

In \cite{Phong_Sturm05}, Phong and Jacob approximate the $C^{1,1}$ geodesic segment
(established in \cite{chen991}) in algebraic manifold via finite dimensional approach.  In light of the preceeding
theorem, it will be nice to approach this relative $C^{1,1}$ geodesic ray via finite dimensional
approaches too.\\

\begin{defi} \label{defi:geodesicinvariant}For every geodesic ray $\rho(t) (t\in [0,\infty)$,  we
can define an invariant as
\begin{equation} \yen(\rho) =  \displaystyle \lim_{t \rightarrow
\infty} \displaystyle \int_M\; {{\p \rho(t)}\over {\p t}}
(\underline{R} - R(\rho(t))) \omega_{\rho(t)}^n.
\label{eq:invariant1}
\end{equation}
\end{defi}
\begin{rem} For a smooth geodesic ray,  the K energy is convex and the above
invariant is well defined.  \\

   In case the geodesic ray arises from a one parameter holomorphic
   transfermation, the integrand in equation \ref{eq:invariant1} is
   just the usual Calabi-Futaki invariant.   This invariant shall be compared to the genearlized Futaki
   invariant defined by Ding and Tian on Fano varieties.

\end{rem}

A natural question is:  If two geodesic rays are parallel to each other, are their
$\yen$ invariants the same?   The answer is partially``yes":
\begin{prop} If one of the geodesic rays has bounded
ambient geometry, then any other geodesic ray parallel to it must have same $\yen$ invariant.
\end{prop}

\begin{defi}\label{def:geodesicdestable}
A geodesic ray $\rho(t): [0,\infty) \rightarrow \cH$ is called
stable (resp; semi-stable) if $\yen(\rho)>0$ (resp: $\geq 0$).   It is called a de-stabalizer for
$\cH$ if $\yen(\rho) < 0\;$ and it is called an effective
de-stablizer if in additional
\[
\displaystyle 
\limsup_{t\rightarrow \infty}\;\; t^2 \cdot \displaystyle \int_M\;  (R(\rho(t))-
\underline{R})^2 \omega_{\rho(t)}^n = 0.
\]

\end{defi}

Following the approach in algebraic case,  we define (cf. \cite{Dona96}):
\begin{defi}\label{def:geodesicstable} A K\"ahler  manifold is called (effectively)  geodesicly stable if
there is no (effective) de-stablizing geodesic ray.  It is called
weakly geodesicly stable if the invariant $\yen$ is always
non-negative for every geodesic ray.
\end{defi}

One of the main theorems is:

\begin{theo}\label{th:calabienergybound} Suppose $\rho(t):[0,\infty) \rightarrow \cH$ is
an effective de-stablized geodesic ray in $\cH$,  then
\[
  \int_M\;(R(\varphi)- \underline{R})^2 \omega_{\varphi}^n  \geq \yen(\rho)^2, \qquad \forall\;\varphi \in \cH.
\]
In fact,  we have \begin{equation} \displaystyle \inf_{\varphi
\in\cH} \int_M\; (R(\varphi)- \underline{R})^2 \omega_\varphi^n
\geq \displaystyle \sup_{\rho} \yen(\rho)^2,
\label{th:calabienergybound}
\end{equation}
where the sup in the right hand side of \ref{th:calabienergybound}
runs over all possible effective de-stablized geodesic rays.
\end{theo}

As a corollary, we have the following important consequence
\begin{cor} If there is a K\"ahler metric of constant scalar
curvature, then it is weakly effectively geodesic stable.
\end{cor}

One can generalized these results to the case of extremal K\"ahler
metric with non-constant scalar curvature.\\

\noindent {\bf
Definition \ref{defi:geodesicinvariant}a} {\it  Suppose $\cX_c$ is
the canonical extremal vector field in $(M,[\omega])$ (cf. eq.
\ref{eq:extremalvectorfield}) and $\theta(\cX_c)$ is defined as
equation \ref{eq:liederivatives}. For every geodesic ray $\rho(t)
(t\in [0,\infty))$, we can define an invariant as
\begin{equation} \tilde\yen(\rho) =  \displaystyle \lim_{t \rightarrow
\infty} \displaystyle \int_M\; {{\p \rho(t)}\over {\p t}}
(\underline{R} - R(\rho(t) - \theta(\cX_c))) \omega_{\rho(t)}^n.
\label{eq:invariant2}
\end{equation}
A geodesic ray $\rho(t): [0,\infty) \rightarrow \cH$ is called
stable (resp; semi-stable) if $\tilde \yen(\rho)>0$(resp: $\geq 0$).  It is called a destabilizer
for $\cH$ if $\tilde\yen(\rho) < 0\;$ and effective de-stablizer
if in addition
\[
\displaystyle \limsup_{t\rightarrow \infty}\;\; t^2 \cdot \int_M\; (R(\rho(t))-
\underline{R}-\theta(\cX_c))^2 \omega_{\rho(t)}^n  =  0.
\]}
With essentially same proof, we have\\

\noindent{\bf Theorem \ref{th:calabienergybound}a} {\it Suppose
$\rho(t):[0,\infty) \rightarrow \cH$ is an effecitvely destabilizing geodesic ray
in $\cH$,  then
\[
  \int_M\;(R(\varphi)- \underline{R} -\theta(\cX_c))^2 \omega_{\varphi}^n  \geq \tilde\yen(\rho)^2,\qquad
  \forall\; \varphi \in \cH.
\]
In fact,  we have \begin{equation} \displaystyle \inf_{\varphi
\in\cH} \int_M\; (R(\varphi)- \underline{R} -\theta(\cX_c))^2
\omega_\varphi^n \geq \displaystyle \sup_{\rho}
\tilde\yen(\rho)^2, \label{th:calabienergybounda}
\end{equation}
where the sup in the right hand side runs over all possible effectively
destabilizing geodesic rays. Moreover, the underlying manifold is
weakly geodesic stable if there exists an extremal K\"ahler metric
in the K\"ahler class.}

\subsection{Proof of Theorem \ref{th:geodesicDegenerate}}
In this subsection, we will give a proof of the existence of a geodesic
ray when  the initial geodesic ray has bounded ambient geometry.  One of the main
challenges  here has been searching for the right condition for the existence of a parallel geodesic ray with regularity beyond the
$L^2$ topology on the K\"ahler potential. 
Following the main steps in \cite{chen991} under current circumstance: for any given
K\"ahler potential $\varphi_0$, we can pick a
sequence of K\"ahler metrics $\rho(t_i) (i \in\NN) $ along the given geodesic ray, 
and connects $\varphi_0$
to $\rho(t_i) (i\in \NN)$ via the unique
$C^{1,1}$ geodesic segment established in \cite{chen991}.   This way, we obtain a sequence of
$C^{1,1}$ geodesic ray and hope to take a limit as $t_i\rightarrow \infty.\;$ The main difficult is to
obtain some uniform $C^{1,1}$ bound which allows us to take a
limit as $t_i\rightarrow \infty.\;$  However, such an approach runs into a serious problem as we
shall explain now: first, there is no absolute $C^0$ estimate
which is crucial to the Yau's calculation of the second
derivatives.  Secondly,  when we do the blowing up estimate, the compactness of
the underlying K\"ahler manifold becomes crucial.  Thirdly, in
deriving boundary estimate as in \cite{chen991}, we need the
assumption that the restriction of K\"aher metric in $\p \Sigma
\times M$ has a uniform positive lower bound with respect to some
fixed metric. This is clearly not available since the sequence of
metrics along a geodesic ray are expected to either diverge or
converge to a metric in different complex structure.  In a
typical senario, this sequence of metrics will degenerate along 
generic points in the K\"ahler manifold and will blowup along some
divisor.  To overcome this difficulty, we use this bounded ambient metric 
from the initial geodesic ray to obtain some control of $C^0$ bound
on the modified potentials.  In order to derive a $C^2$ estimate
in terms of this weak $C^0$ estimate, we need to exploit the structure
of degenerated Monge-Ampere equation more closely. In particular, if the
modified potential doesn't have a uniform $C^0 $ bound, we need
to re-design the blowing up procedure in \cite{chen991} to obtain a
growth control in the $C^{1,1}$ bound on the modified potential.  We believe that
such a technique may be applicable to some other interesting cases.
%slightly to address the issue that there is no absolute $C^0$
%estimate on the K\"ahler potential. 
%Moreover,  we need to understand the compactness of a family of
%manifolds $(M, \omega_{\rho(t)}).\;$ Such a compactness is usually
%not known unless this geodesic ray is corresponding to some
%special degeneration. In general, we need to
%come up some other ideas to overcome this difficulties.  As stated before,  we regard
%this as a first step towards a general existence theorem on geodesic ray in
%subsequent paper. \\

\subsubsection{Setup of problem}
Let us first set up some notations. Let $T \gg 1 $ be a generic
large positive number. Let  $\Sigma_T = [0,T]\times S^1.\;$ In the
$(n+1)$ dimensional K\"ahler manifold $\Sigma_T\times M$, we want
to solve the Dirichelet problem for  HCMA equation \ref{eq:hcma0}
where the boundary data is invariant in the circle $S^1$
direction. As in \cite{chen991},  for any $T$ and for any smooth
boundary data, we can obtain a unique $C^{1,1}$ solution $\phi(t)$
such that it solves HCMA equation \ref{eq:hcma0}.  In other words,
we have
\begin{equation}
\left( \pi_2^* \omega_0 + \sqrt{-1} \p \b \p \phi\right)^{n+1} =
0, \label{eq:hcma00} \end{equation}
where
\begin{equation}\phi(0) =\varphi_0, \qquad {\rm and}\qquad \phi(T) = \rho(T).
\label{eq:hcma00a}
\end{equation}
Here we assume that $\rho:[0,\infty) \rightarrow \infty$ is a
smooth geodesic ray.  Obviously, the $C^{1,1}$ estimate depends on
$T$ and may blow up as $T \rightarrow \infty.\;$ In fact,  this
$C^{1,1}$ estimate must blow up if it represents a geodesic ray.
This creates a serious problem for the existence of geodesic rays.
Our strategy is the following: let $\{n_k\in \NN\}$ be a sequence
of numbers that approach $\infty.\; $ Suppose that $\{\phi_k, k\in
\NN\}$ solves the Dirichelt boundary value problem of equation
\ref{eq:hcma0} in $\Sigma_{n_k}\times M$ with boundary data
\begin{equation}\phi_k(0) =\varphi_0, \qquad {\rm and}\qquad \phi_k(n_k) = \rho(n_k).
\label{eq:hcma00k}
\end{equation}

\begin{lem} \label{lem:c0estimate}For any smooth geodesic ray $\rho(t) (t\in
[0,\infty))$ and for any initial metric $\varphi_0 \in \cH$,
there exists a uniform constant $C$ such that for any $T\in
(0,\infty)$, there exists a unique $C^{1,1}$ geodesic  $\phi_T(t)
(t\in [0, T])$ which connects $\varphi_0$ to $\rho(T)$ such
that
\begin{equation}
- C \leq \phi_T(t,x) - \rho(t,x)  \leq C.
\label{eq:c0estimate}\end{equation}
\end{lem}

To obtain uniform $C^{1,1}$ bound in some fashion, we need to choose some
appropriate background K\"ahler metric first.
Let $h$ be the ambient metric  with bounded ambient geometry. Suppose that this initial geodesic ray $\rho(t) (t\in [0,\infty))$ is tamed by $h.\;$
Suppose that its K\"ahler form $\tilde{\omega}$ is given
by
\begin{equation}
  \pi_2^* \omega_0 + \displaystyle \sum_{i,j=1}^n\; {{\p^2 \b \rho}\over
  {\p w^i \p w^{\b j}}} + 2 Re\left(\displaystyle \sum_{i=1}^n\; {{\p^2 \b \rho}\over
  {\p w^i \p \b z}} d\, w^i d\, w^{\b j} \right) +  {{\p^2  \b \rho}\over {\p z \p \b z}}
  d\, z\; d\, \b z.
\label{eq:wrappedproductmetric} \end{equation} Here $z = t +
\sqrt{-1} \theta.\;$ In other words
\[
\tilde{\omega} = \pi_2^*\omega_0 + \sqrt{-1} \p \b \p \b \rho .
\]

The Dirichelet boundary value problem eq. \ref{eq:hcma00} and  \ref{eq:hcma00a} can be re-written as
  a Drichelet problem on
  $\Sigma_T\times M$ such that
\begin{eqnarray}
\det\left( h_{\alpha \b \beta} + {{\p^2 (\phi- \b \rho)}\over {\p w^{\alpha}\p w^{\b \beta}}}\right)_{(n+1)\times (n+1)}
= 0, \label{eq:hcma1}
\end{eqnarray}
with boundary condition
\begin{equation}
  \phi\mid_{\{0\}\times S^1 \times M} =  \varphi_0, \;\;{\rm
  and}\;\;  \phi\mid_{\{T\}\times S^1 \times M} =  \rho(T).
  \label{eq:hcmaT}
\end{equation}
%{\bf The boundary condition should be}
%\[ \phi\mid_{\{T\}\times S^1 \times M} =  \rho(T). \]

%Denote the solution of this problem as $\phi_T(t,x).\; $
Set
\begin{equation}
\tilde \psi_T(t,x)= \phi_T(t,x) -\b \rho(t,x). 
\label{eq:modifiedpotential}
\end{equation}
For a sequence of points $t_i \rightarrow \infty$ we have
\[
\tilde \psi_{t_i}(t_i, x) = \phi_{t_i}(t_i, x) - \b \rho(t_i,x) = \rho_{t_i}(t_i, x) - \b \rho(t_i, x).
\]
%and \begin{equation} \tilde \psi_T(t,x) = \psi_T(x)  +  k(t). \label{eq:modifiedpotential1} \end{equation}
For simplicity,
we drop the dependency on $T.\;$ Thus, the modified potential
$\tilde \psi_T(t,x) $ has  uniform $C^0$
bound. % but $\tilde\psi(t,x)$ does!\\

As in
\cite{chen991}, we want to use the method of continuity. So we set up the
problem as\\
\begin{eqnarray}
\det\left( h_{\alpha \b \beta} + {{\p^2 \tilde \psi}\over {\p w^{\alpha}\p
w^{\b \beta}}}\right)_{(n+1)\times (n+1)} = \epsilon \det
\left(g_{i\b j}\right)_{n\times n}, \label{eq:hcma2}
\end{eqnarray}

with boundary condition
\begin{equation}
 \tilde  \psi\mid_{\{0\}\times S^1 \times M} = \varphi_0 -\rho(0) , \;\;{\rm
  and}\;\;  \tilde \psi\mid_{\{T\}\times S^1 \times M} =\rho(T) - \b \rho(T).
  \label{eq:modifiedpotentialA}
\end{equation}
For any $T$ fixed, this Drichelet boundary value has a unique
$C^{1,1}$ solution as in \cite{chen991}. The challenge at hand is
how to obtain  a $C^{1,1}$ estimate when $T$ runs over an increasing sequence
of times $\{n_k\in \NN\}$ and $\epsilon \rightarrow 0.\;$

\subsubsection{The $C^{1,1}$ estimates for the HCMA equation in unbounded domains}

In this subsection, we want to solve equation \ref{eq:hcma2} for any large $T>0.\;$ We follow Yau's estimate in \cite{Yau78} and we want to set up
some notations first.
Put $\omega_{\tilde{\rho}(t)} = \sqrt{-1}
h_{\alpha \bar \beta} d\,w^{\alpha} \otimes w^{\bar \beta}\;$ and $
\omega_{\varphi(t)} = \sqrt{-1} g'_{\alpha \bar \beta} d\,w^{\alpha}
\otimes d\,w^{\bar \beta} \;$ where \[ g'_{\alpha\bar \beta} =
h_{\alpha \bar \beta} + {{\p^2 \left(\varphi(t) -\b {\rho}(t) \right)
}\over {\p w^{\alpha} \p w^{\bar \beta}}}.
\] Then
\[
  \triangle' =  \displaystyle \sum_{\alpha,\beta=1}^n\; {g'}^{\alpha \bar \beta} {{\partial^2 }\over
  {\partial w^{\alpha} \partial w^{\bar \beta}}},\qquad
   \tilde{\triangle} = \displaystyle \sum_{\alpha,\beta=1}^n\; h^{\alpha \bar \beta} {{\partial^2\;}\over
  {\partial w^{\alpha} \partial w^{\bar \beta}}}.
\]
Before stating the crucial Lemma of this subsection, we need to explain a little bit the relationship
between geodesic ray and its bounded ambient geometry.  By definition of the initial geodesci ray 
tamed by an ambient metric $h$ with bounded ambient geometry, there exists a sequence
of time $t_i \rightarrow \infty$ such that
%\[
%\rho(t_i) = \b \rho(t_i), \qquad \forall \; i=1,2,\cdots.
%\]
%Moreover, the following hold
\begin{enumerate}
\item $|n+1 +  \triangle_h (\rho - \b \rho)|  \leq C;$
\item $  |{{\p (\rho -\b \rho)}\over {\p t}} |_h \leq C;$
\item The vector $|{\p\over {\p t}}|_h$ has uniform upper bound. 
\end{enumerate}
Note that under the first two conditions, $\rho-\b \rho$ is not necessary bounded.  However, it is sufficient to show that the oscillation of $|\rho-\b \rho|$ is controlled by the distance (by ambient metric $h$).  
We first prove the following lemma:
%\begin{lem} \label{lem:yau} If the $C^0$ norm of  $\mid \tilde \psi\mid$ is uniformly bounded (independent of time $t$), then there
%exists a uniform constant $C$ such that
%\[ \triangle'u \geq c_0 u^{n\over {n-1}} - c_1 u - c_2 \]
%where \[
%u = e^{-C\tilde \psi} (n + 1 + \tilde{\triangle} \psi).
%\]
 % All constants are
%independent of $T$ if the K\"ahler metric $h$ in $\Sigma_T \times M$
% has  a uniform lower bound on its bisectional
%curvature.

%\end{lem}
\begin{lem} \cite{chen991}\label{cor:2ndderivativeestimate} There exists a constant $C$ which depends only on the ambient metric $ h$ (independent of $T$) such that
\[
e^{-\lambda (\rho -\b \rho)}
( n + 1+ \tilde{\triangle} \tilde \psi(t) ) \leq \displaystyle \max_{t=0, t=T} \; e^{-\lambda (\rho -\b \rho)}
( n + 1+ \tilde{\triangle} \tilde \psi(t) ). \]
%Either\footnote{\label{growthcondtion on h}Here we assume $|{\p \over {\p t}}|_h$ has a positive uniform lower bound.
%The inequality hold for any $t< T$ such that minimal distance in $h$ from $\{t\} \times S^1 \times M$ to
%$\{T\}\times S^1\times M$ is finite (say less than $1$).  The important matter is the oscillation
%of $\rho-\b \rho$ must be controlled. }  
%\begin{equation}
%\displaystyle \max_{ [T-1, T]\times S^1 \times M} \; ( n + 1+ \tilde{\triangle} \tilde \psi(t) )
% \leq C \displaystyle \max_{\{T\}\times S^1  \times M} \;( n + 1 + \tilde{\triangle} \tilde \psi(t) ).
% \label{eq:2ndderivativecontrol1} \end{equation}
%or \begin{equation} \displaystyle \max_{[0,1]\times S^1 \times M} \;
 %( n + 1+ \tilde{\triangle} \tilde \psi(t) )\leq C \displaystyle \max_{\{0\}\times S^1 \times M} \;
% ( n + 1 + \tilde{\triangle} \tilde \psi(t) ).\label{eq:2ndderivativecontrol2}
%\end{equation}
%holds.   The choice depends where $e^{-\lambda (\rho -\b \rho)}
%( n + 1+ \tilde{\triangle} \tilde \psi(t) ) $ attains the maximal value.
\end{lem}

\begin{proof} We want to use the maximum principle in this proof.
Let us first calculate $ \triangle'\left( n + 1 +
\tilde{\triangle} (\varphi-\b \rho) \right).\;$

Let us choose a coordinate so that at a fixed point both
$\omega_{\tilde{\rho}(t)} =\sqrt{-1} h_{\alpha \bar \beta}
d\,w^{\alpha} \otimes d\,w^{\bar \beta}\; $  and the complex Hessian
of $\varphi(t)-{\rho}(t)$ are in diagonal forms. In
particular, we assume that $ h_{i \bar j} = \delta_{i \bar j} $ and
$ \left(\varphi(t) -{\rho}(t)\right)_{i \bar j} = \delta_{i
\bar j} \left(\varphi(t) -{\rho}(t)\right)_{i \bar i}.\;$ Thus
\[
{g'}^{ i \bar s} = {{\delta_{i \bar s}} \over {1 + \left(\varphi(t)
-{\rho}(t)\right)_{i \bar i} }}.\]

For convenience, put
\[
F = \ln \epsilon + \log \det \left(h_{i\b j}\right).
\]

Then our equation reduces to
\[
\log \det\left(h_{i \bar j} + {{\p^2 \left(\varphi - \b {\rho}\right)} \over {\p w_{i} \p w_{\bar j}}} \right)
= F + \log \det (h_{i \bar j}).
\]
For convenience, set
\[
   \tilde  \psi(t) = \varphi(t) -\b { \rho}(t)
\]
in this proof. Note that $\mid \tilde\psi(t)\mid$ is  uniformly bounded.
We first follow the standard calculation of $C^2$ estimates in \cite{Yau78}.
Differentiate both sides with respect to ${\p \over {\p w_k}}$
$$
(g')^{i \b j} \biggl( {{\p h_{i \b j}} \over {\p w_k}} +{{\p ^3 \tilde
\psi(t)} \over { \p w_i \p \b w_j \p w_k}} \biggr) - h^{i \b j} {{\p
h_{i \b j}} \over {\p w_k}} = {{\p F} \over {\p w_k}},
$$
and differentiating again with respect to ${\p \over {\p \b w_l}}$
yields
$$
(g')^{i \b j} \biggl( {{\p ^2 h_{i \b j}} \over {\p w_k \p \b w_l}}
+ {{\p ^4 \tilde \psi(t)} \over { \p w_i \p \b w_j \p w_k \p \b w_l}}
\biggr) +h^{t \b j} h^{i \b s} {{\p h_{t \b s}} \over {\p \b w_l}}
{{\p h_{i \b j}} \over {\p w_k}} -h^{i \b j} {{\p ^2 h_{i \b j}}
\over {\p w_k \p \b w_l}}
$$
$$
- (g')^{t \b j}(g')^{i \b s} \biggl( {{\p h_{t \b s}} \over {\p \b
w_l}} +{{\p ^3 \tilde \psi(t)} \over { \p w_t \p \b w_s \p \b w_l}} \biggr)
\biggl( {{\p h_{i \b j}} \over {\p w_k}} +{{\p ^3 \tilde \psi(t)} \over {
\p w_i \p \b w_j \p w_k}} \biggr)
= {{\p ^2 F} \over {\p w_k \p \b w_l}}.
$$
Assume that we have normal coordinates at the given point, i.e.,
$h_{i \b j} = \delta _{ij}$ and the first order derivatives of $g$
vanish. Now taking the trace of both sides results in
\begin{align*}
\Tilde{\Delta} F & = h^{k \b l} (g')^{i \b j}\biggl( {{\p ^2 h_{i \b
j}} \over {\p w_k \p \b w_l}} +{{\p ^4 \tilde \psi(t)} \over { \p w_i \p \b
w_j \p w_k \p \b w_l}}
\biggr) \\ & \qquad -h^{k \b l}
(g')^{t \b j}(g')^{i \b s}{{\p ^3 \tilde \psi(t)} \over { \p w_t \p \b w_s
\p \b w_l}} {{\p ^3 \tilde \psi(t)} \over { \p w_i \p \b w_j \p w_k}} -h^{k
\b l} h^{i \b j} {{\p ^2 h_{i \b j}} \over {\p w_k \p \b w_l}} .
\end{align*}
On the other hand, we also have
\begin{align*}
\Delta ' (\Tilde{\Delta} \tilde \psi(t)) &= (g')^{k \b l} {{\p ^2 } \over
{\p w_k \p \b w_l}}
         \biggl( h^{i \b j} {{\p ^2 \tilde \psi(t)} \over {\p w_i \p \b w_j}}\biggr) \\
      &=(g')^{k \b l} h^{i \b j}  {{\p ^4 \tilde \psi(t)} \over { \p w_i \p \b w_j
          \p w_k \p \b w_l}} +(g')^{k \b l} {{\p ^2 h^{i \b j} } \over {\p w_k
         \p \b w_l}} {{\p ^2 \tilde \psi(t)} \over {\p w_i \p \b w_j}},
\end{align*}
and we will substitute ${{\p ^4 \tilde \psi(t)} \over {\p w_i \p \b w_j\p
w_k \p \b w_l}} $ in $\Delta ' (\Tilde{\Delta} \tilde \psi(t))$ so that the
above reads
\begin{align*}
\Delta ' (\Tilde{\Delta}  \tilde \psi(t)) &= -h^{k \b l} (g')^{i \b j} {{\p
^2 h_{i \b j}} \over
              {\p w_k \p \b w_l}} +h^{k \b l}
(g')^{t \b j}(g')^{i \b s}{{\p ^3  \tilde \psi(t)} \over { \p w_t \p \b w_s
\p \b w_l}}
{{\p ^3 \tilde \psi(t)} \over { \p w_i \p \b w_j \p w_k}} \\
&\quad +h^{k \b l} h^{i \b j}
{{\p ^2 h_{i \b j}} \over {\p w_k \p \b w_l}} +\Tilde{\Delta} F
+(g')^{k \b l} {{\p ^2 h^{i \b j} } \over {\p w_k
         \p \b w_l}} {{\p ^2 \tilde \psi(t)} \over {\p w_i \p \b w_j}},
\end{align*}
which we can rewrite after substituting ${{\p ^2 h_{i \b j}} \over
{\p w_k \p \b w_l}} = -R_{i \b j k \b l}$ and  ${{\p ^2 h^{i \b j}}
\over {\p w_k \p \b w_l}} = R_{j \b i k \b l}$ as
$$
\begin{array}{lcl}
\Delta ' (\Tilde{\Delta}  \tilde \psi(t)) & = & \Tilde{\Delta} F +h^{k \b
l}(g')^{t \b j}(g')^{i \b s} \tilde \psi(t) _{t \b s l} \tilde \psi(t) _{i \b j k}
\\ & & +(g')^{i \b j} h^{k \bar l} R_{i \b j k \b l} - h^{i \b j}
h^{k \bar l}R_{i \b j k \b l} +(g') ^{k \b l} R_{j \b i k \b l}
\tilde \psi(t) _{i \b j}.
\end{array}
$$
Restrict to  the coordinates we chose in the beginning so that both
$g$ and $\tilde \psi(t)$ are in diagonal form. The above transforms to
$$
\begin{array}{lcl}
\Delta ' (\Tilde{\Delta} \tilde \psi(t)) & = & {1 \over {1 +\tilde \psi(t) _{i \b
i}}} {1 \over {1 + \tilde \psi(t) _{j \b j}}} \tilde \psi(t) _{i \b j k} \tilde  \psi(t)
_{\b i j \b k} + \Tilde{\Delta} F \\ & & \qquad \qquad +R_{i \b i k
\b k} (-1+ {1 \over {1+\tilde \psi(t) _{i \b i}}}+{{\tilde \psi(t) _{i \b i}}
\over {1+\tilde \psi(t) _{k \b k}}}).
\end{array}
$$
Set now $C= \inf _{i \ne k} R_{i \b i k \b k}$ and observe that
\begin{align*}
R_{i \b i k \b k} (-1+ {1 \over {1+\tilde \psi(t) _{i \b i}}}+{{\tilde \psi(t) _{i \b i}}
\over {1+\tilde \psi(t) _{k \b k}}})
  &= {1\over 2} { R_{i \b i k \b k}}{{(\tilde \psi(t) _{k \b k} - \tilde\psi(t) _{i \b i})^2}
     \over {(1 +\tilde \psi(t)_{i \b i})(1+\psi(t) _{k \b k})}} \\
  &\ge {C \over 2} {{(1+\tilde \psi(t) _{k \b k} -1-\psi(t) _{i \b i})^2} \over {(1 +\tilde \psi(t)
   _{i \b i})(1+\tilde \psi(t) _{k \b k})}} \\
  &= C\Bl {{1 +\tilde \psi(t) _{i \b i}} \over {1+\tilde \psi(t) _{k \b k}}} -1
  \Br,
\end{align*}
which yields
$$
\begin{array}{lcl}
\Delta ' (\Tilde{\Delta} \tilde \psi(t)) & \ge & {1 \over {(1 + \tilde \psi(t) _{i
\b i})(1 +\tilde \psi(t) _{j \b j})}} \tilde \psi(t) _{i \b j k} \tilde \psi(t) _{\b i j
\b k} + \Tilde{\Delta}
F \\
& & \qquad \qquad +C \Bl (n + 1+ \Tilde{\Delta} \tilde \psi(t)) \sum_i {1
\over {1 + \tilde \psi(t) _{i \b i}}}-1 \Br.
\end{array}
$$
We need to apply one more trick to obtain the requested estimates.
Namely,
\[
\begin{array} {lcl} & &
\Delta ' (e^{-\l \tilde \psi(t)} (n + 1 + \Tilde{\Delta}\tilde \psi(t)) )\\
  &= &  e^{-\l \tilde \psi(t)} \Delta ' (\Tilde{\Delta} \psi(t)) +2\nabla ' e^{-\l\tilde \psi(t)} \nabla '
   (n + \Tilde{\Delta} \psi(t)) \\
     & & \qquad +\Delta ' ( e^{-\l \tilde \psi(t)}) (n + 1 + \Tilde{\Delta} \tilde \psi(t)) \\
  &= & e^{-\l \tilde \psi(t)} \Delta ' (\Tilde{\Delta}  \tilde \psi(t)) -\l  e^{-\l \tilde \psi(t)} (g')^{i \b i}
   \tilde \psi(t) _i (\Tilde{\Delta}\tilde \psi(t) )_{\b i}
     \\ & & \qquad -\l  e^{-\l \tilde \psi(t)} (g')^{i \b i} \tilde \psi(t) _{\b i}
     (\Tilde{\Delta}\tilde \psi(t)) _i \\
  & & \quad - \l  e^{-\l \tilde \psi(t)} \Delta' \tilde \psi(t) (n + 1 + \Tilde{\Delta}\tilde \psi(t))
  \\ & & +\l ^2 e^{-\l \tilde \psi(t)}
    (g')^{i \b i} \tilde \psi(t) _{i} \tilde \psi(t)_{ \b i}  (n + 1 +\Tilde{\Delta}\tilde \psi(t)) \\
  &\ge &  e^{-\l \tilde \psi(t)} \Delta ' (\Tilde{\Delta}\tilde \psi(t))\\
  & &  -e^{-\l \tilde \psi(t)}(g')^{i \b i}
    (n + \Tilde{\Delta} \tilde \psi(t))^{-1}(\Tilde{\Delta}  \tilde \psi(t) )_i (\Tilde{\Delta} \tilde \psi(t) ) _{\b i} \\
  & &\quad -\l e^{-\l \tilde \psi(t)}\Delta ' \tilde \psi(t) (n + 1 + \Tilde{\Delta}\tilde \psi(t)) ,
\end{array}
\]
which follows from the Schwarz Lemma applied to the middle two
terms. We will write out one term here; the other goes in an
analogous way.
\begin{align*}
&(\l e^{-{\l \over 2} \tilde \psi(t)} \tilde \psi(t)_i
(n+\Tilde{\Delta} \tilde \psi(t) )^{{1 \over 2}} )
( e^{-{\l \over 2} \tilde \psi(t)} (\Tilde{\Delta}\tilde \psi(t)) _{\b i}(n + 1 +\Tilde{\Delta}\tilde\psi(t) )^{-{1 \over 2}})\\
\le & {1 \over 2} (\l ^2 e^{-\l \tilde \psi(t)} \tilde \psi(t) _i
\tilde \psi(t) _{\b i} (n + 1 +\Tilde{\Delta}
\tilde \psi(t) ) \\
& +e^{-\l \tilde \psi(t)}(\Tilde{\Delta}\tilde\psi(t)) _{\b i}
(\Tilde{\Delta} \tilde \psi(t))_i
   (n+\Tilde{\Delta}\tilde \psi(t) )^{-1}).
\end{align*}
Consider now the following
\begin{align*}
&-(n + 1 +\Tilde{\Delta} \tilde \psi(t) )^{-1} {1 \over {1+ \tilde \psi(t) _{i \b
i}}} (\Tilde{\Delta} \tilde \psi(t)) _i (\Tilde{\Delta}
\tilde \psi(t) )_{\b i} +\Delta ' \Tilde{\Delta} \tilde \psi(t) \ge \\
& -(n + 1 +\Tilde{\Delta}\tilde \psi(t) )^{-1} {1 \over {1+\tilde \psi(t) _{i
\b i}}} | \tilde \psi(t) _{k \b k i}|^2
+\Tilde{\Delta} F \\
&+ {1 \over {1+\tilde \psi(t) _{i \b i}}} {1 \over {1+ \tilde \psi(t) _{k \b k}}}
\tilde \psi(t) _{k \b i \b j}  \tilde \psi(t) _{i \b k j} +C(n + 1
+\Tilde{\Delta} \tilde \psi(t) ){1 \over {1+\tilde \psi(t) _{i \b i}}}.
\end{align*}
On the other hand, using the Schwarz inequality, we have
\begin{align*}
& (n+\Tilde{\Delta} \tilde \psi(t) )^{-1}{1 \over {1+\tilde \psi(t) _{i \b i}}} |\tilde \psi(t) _{k \b k i}|^2 \\
& \quad \quad =(n + 1 +\Tilde{\Delta} \tilde \psi(t) )^{-1}{1 \over
{1+\tilde \psi(t) _{i \b i}}} \Biggl|{{ \tilde \psi(t) _{k \b k i} } \over
{(1+ \tilde \psi(t) _{k \b k} )^{1 \over 2}}}
(1+\tilde \psi(t) _{k \b k} )^{1 \over 2} \Biggr| ^2 \\
&\quad \quad \le (n + 1 +\Tilde{\Delta} \psi(t) )^{-1} \Bl {1
\over {1+\tilde \psi(t) _{i \b i}}}
{1 \over {1+ \tilde \psi(t) _{k \b k}}}\tilde\psi(t) _{k \b k i} \tilde \psi(t) _{\b k k \b i} \Br
\Bl 1+ \tilde \psi(t) _{l \b l} \Br \\
&\quad \quad ={1 \over {1+ \tilde \psi(t) _{i \b i}}}
{1 \over {1+\tilde\psi(t) _{k \b k}}} \tilde \psi(t) _{k \b k i}\tilde \psi(t) _{\b k k \b i} \\
&\quad \quad ={1 \over {1+ \tilde \psi(t) _{i \b i}}}{1 \over {1+\tilde \psi(t) _{k
\b k}}}
\psi(t) _{i \b k k } \tilde \psi(t) _{k \b i \b k} \\
&\quad \quad \le {1 \over {1+\tilde \psi(t) _{i \b i}}}{1 \over {1+\tilde \psi(t)
_{k \b k}}}
\tilde\psi(t) _{i \b k j }\tilde \psi(t) _{k \b i \b j},
\end{align*}
so that we get
$$
\begin{array}{l}
-(n+\Tilde{\Delta}\tilde \psi(t) )^{-1} {1 \over {1+\tilde\psi(t) _{i \b i}}}
(\Tilde{\Delta} \tilde\psi(t)) _i (\Tilde{\Delta} \tilde\psi(t) )_{\b i}
+\Delta ' \Tilde{\Delta} \tilde\psi(t) \\
\ge \Tilde{\Delta} F + C(n + 1 +\Tilde{\Delta} \tilde\psi(t) ) {1 \over
{1 + \tilde\psi(t) _{i \b i}}}.
\end{array}
$$
Putting all these together, we obtain

\begin{eqnarray}
\triangle' \left(e^{-\lambda \tilde \psi(t) } (n + \tilde{\triangle} \tilde \psi(t)
)\right) \qquad \qquad \nonumber \\
\geq e^{-\lambda \tilde \psi(t) }
\left( \tilde{\triangle} F + C  (n + 1 + \tilde{\triangle} \tilde\psi(t))
  \displaystyle \sum_{i=1}^n\;{1 \over {1 + \tilde \psi(t)_{i \bar i} }}
  \right) \nonumber \\
   \qquad \qquad - \lambda\; e^{-\lambda  \tilde \psi(t)
  } \triangle' \tilde \psi(t) \; (n + 1 + \tilde{\triangle} \tilde \psi). \qquad
  \label{eq:estimate1}
\end{eqnarray}

Consider
\[
\begin{array}{lcl} \tilde{\triangle} F & = & h^{\alpha \b \beta} {{\p^2 \log \det (h_{i\b j})}\over {\p w^\alpha \p w^{\b \beta}}}
= - R(\rho(t)).
\end{array}
\]

Plugging this into the inequality (\ref{eq:estimate1}), we obtain
\[ \begin{array}
{l}\triangle'  \left( e^{-\lambda \tilde \psi(t) } ( n +
\tilde{\triangle} \psi(t) )
\right) \\
\geq e^{-\lambda \tilde\psi(t) }
\left(  C  (n + 1 +\tilde{\triangle} \psi(t))
  \displaystyle \sum_{i=1}^n\;{1 \over {1 + \tilde\psi(t)_{i \bar i} }}
  \right) \nonumber \\
   \qquad \qquad - \lambda\; e^{-\lambda \tilde\psi(t)
  } \triangle' \tilde \psi(t) \; (n + 1 + \tilde{\triangle}  \tilde\psi(t))  - R(\rho(t))\; e^{-\lambda \tilde\psi(t) }.
\end{array}\]

Now
\[ \begin{array} {lcl}\triangle' \tilde \psi(t) & = &  \triangle' \tilde  \psi(t)
= tr_{g'} (\tilde \omega + i \p \b \p \tilde\psi - (\tilde \omega )
\\& = & n + 1  - 
tr_{g'} h.\end{array}\]
%Here
%\[\tilde h = \pi_2^* \omega_0 + i \p \b \p \b \rho > \delta h.\]
%Then, \[
%\triangle' \tilde \psi(t) \leq n+1 - \delta   \displaystyle \sum_{i=1}^{n+1}\;{1 \over {1 + \psi(t)_{i \bar i} }}.\]
Plugging this into the
above inequality, we obtain
\[\begin{array}{l} \triangle'   \left( e^{-\lambda \tilde \psi(t) }
( n + \tilde{\triangle} \tilde \psi(t) ) \right) \\ \geq e^{-\lambda \tilde
\psi(t) }
\left(  (C + \lambda \delta)  (n +1+ \tilde{\triangle} \tilde \psi(t))
  \displaystyle \sum_{i=1}^{n+1}\;{1 \over {1 + \psi(t)_{i \bar i} }}
  \right) \nonumber \\\qquad - |\lambda|
c_3 \;e^{-\lambda  \tilde\psi(t) } \; ( n + 1 + \tilde{\triangle} \tilde
\psi(t) ) - R(\rho(t))\; e^{-\lambda \tilde \psi(t) }.
\end{array}
\]
Let $\lambda \delta = - C + 1,$ we then have
\[ \begin{array}{l}\triangle'  \left( e^{-\lambda \tilde \psi(t) }
( n + \tilde{\triangle} \tilde \psi(t) ) \right) \\ \geq e^{-\lambda \tilde
\psi(t) }
\left( (n + \tilde{\triangle}  \tilde \psi(t))
  \displaystyle \sum_{i=1}^{n+1}\;{1 \over {1 + \tilde \psi(t)_{i \bar i} }}
  \right) \nonumber \\\qquad  \qquad - c_5\;e^{-\lambda \tilde \psi(t) } \; ( n
  + 1+
\tilde{\triangle} \tilde \psi(t) ) - c_2\; e^{-\lambda \tilde\psi  }.
\end{array}\]
Here $c_5$ is a uniform constant.  \\

\noindent {\bf Claim}: the maximum value of
$e^{-\lambda \tilde \psi(t) }
( n + 1 + \tilde{\triangle} \tilde \psi(t) ) $ must occur in  $\p \Sigma_T \times M.\;$\\

Otherwise, if the maximum occur in the interior, we have
\[
 e^{-\lambda \tilde
\psi(t) }
\left( (n + 1+ \tilde{\triangle}  \tilde \psi(t))
  \displaystyle \sum_{i=1}^{n+1}\;{1 \over {1 + \tilde \psi(t)_{i \bar i} }}
  \right) \nonumber - c_5\;e^{-\lambda \tilde \psi(t) } \; ( n
  + 1+
\tilde{\triangle} \tilde \psi(t) ) - c_2\; e^{-\lambda \tilde\psi  } \leq 0.
\]

However, $$\displaystyle \sum_{i=1}^{n+1}\;{1 \over {1 + \tilde \psi(t)_{i \bar i} }}\rightarrow \infty\;$$
as $\epsilon \rightarrow 0.\;$ This leads to a contradiction when $T$ is finite since
\[
\tilde \psi = \rho -\b \rho + \varphi_T -\rho
\]   
is uniformly bounded in $\Sigma_T \times M.$ Thus,

\[
e^{-\lambda \tilde \psi(t) }
( n + 1+  \tilde{\triangle} \tilde \psi(t) ) \leq \displaystyle \max_{t=0, t=T} \; e^{-\lambda \tilde \psi(t) }
( n + 1+ \tilde{\triangle} \tilde \psi(t) ).
\]
In other words, 
\[
e^{-\lambda (\rho -\b \rho)}
( n + 1+ \tilde{\triangle} \tilde \psi(t) ) \leq \displaystyle \max_{t=0, t=T} \; e^{-\lambda (\rho -\b \rho)}
( n + 1+ \tilde{\triangle} \tilde \psi(t) ).
\]
Note that 
\[
|\tilde \triangle (\rho -\b \rho)|  + | {{\p (\rho-\b \rho)}\over {\p t}}| \leq C. 
\]
This implies Lemma 3.17 (cc. \cite{chen991}) 
\end{proof}

As in \cite{chen991}, we have
\begin{theo}\cite{chen991} \label{th:1dcontrolled2d}If $\psi$ is a solution of equation (\ref{eq:hcma2}) at $0<\epsilon <1,$
then  there exists a constant $C$ which depends only on
$(\Sigma_T\times M,h)$ such that if  $e^{-\lambda (\rho -\b \rho)}
( n + 1+ \tilde{\triangle} \tilde \psi(t) ) $ attains the maximal value at $t=T$, then for any $\{t\}\times S^1 \times  M$ which has $h-$ distance to $\{T\}\times S^1 \times M$ less than $1$, we have
\begin{equation}
\displaystyle \max_{\{t\}\times S^1 \times M}\; (n+1 + \tilde \triangle \tilde
\psi) \leq C \displaystyle \max_{[T-\mu, T]\times S^1 \times M} \;(|\nabla \tilde
\psi|_h^2+1),\label{eq:1dderivativecontrol2d1}
\end{equation}
hold  for any $\mu > 0$ where the $h$ distance from $\{T-\mu\}\times S^1 \times M$ to $\{T\}\times M$
is small ($\ll 1$).  On the other hand,  if
$e^{-\lambda (\rho -\b \rho)}
( n + 1+ \tilde{\triangle} \tilde \psi(t) ) $ attains the maximal value at $t=0,$ then
\begin{equation}
\displaystyle \max_{[0,1]\times S^1 \times M}\; (n+1 + \tilde \triangle \tilde
\psi) \leq C \displaystyle \max_{[0,\mu]\times S^1 \times M} \;(|\nabla \tilde
\psi|_h^2+1).\label{eq:1dderivativecontrol2d}
\end{equation}
\end{theo}

%For simplicity, here we assume that $|{\p \over {\p t}}|_h$ has a uniform lower bound.
For simplicity, denote the $h$ distance between two hypersurfaces $\{t-1\}\times S^1 \times M$
and $\{t_2\}\times S^1 \times M$ as $d_h(t_1,t_2).\;$
Following a blowing up argument in \cite{chen991}, we can prove
that  there is a uniform $C^{1,1}$ estimate for $t\in [0,T]\;$ or $t \in [0,1]\;$, depending on where
\[
e^{-\lambda \tilde \psi(t) }
( n +1 + \tilde{\triangle} \tilde \psi(t) )
\]
realizes its maximum.  For simplicity, let us assume that 
\[
\e^{-\lambda \tilde \psi(t) }
( n + 1+  \tilde{\triangle} \tilde \psi(t) )
\]
obtain maximum at $\{T\} \times S^1 \times M.\;$ Thus,  we have 
\[\displaystyle \max_{\{t\}\times S^1 \times M}\; (n+1 + \tilde \triangle \tilde
\psi) \leq C \displaystyle \max_{ [T-\mu, T]\times S^1 \times M} \;(|\nabla \tilde
\psi|_h^2+1), \qquad \forall \mu > 0 \]
for any $t$ where the  $d_h(t, T) \leq 1\;$ and $d_h(T-\mu, T) \ll 1.$   Unlike in \cite{chen991}, here we need
to blowup at the maximam point of
\[
   |\nabla \tilde \psi|_h \cdot ({1\over 2}- d_h(t,T)), \qquad \forall t \in [{T\over 2}, T].
\]
Note that we don't assume that $\tilde \psi$ has a uniform $C^0$ bound.  To bypass this difficulty,
we note that 
\[
\tilde \psi = \varphi - \rho + \rho -\b \rho.
\]
By assumption, the first and second derivatives of the two functions $\tilde \psi$
and $\varphi -\rho$ are equivalent.  Therefore, we really blowup
at the maximum of 
\[
   |\nabla  (\varphi -\rho) |_h \cdot t ({1\over 2}- d_h(t,T)), \qquad \forall t \in [{T\over 2}, T].
\]
As in \cite{chen991}, we can prove that $ |\nabla (\varphi-\rho)|_h $ is uniformly bounded
for $t \in [T-{1\over 2}, T].\;$ Consequently, this implies uniform control on  $ |\nabla \tilde \psi|_h.\;$
The crucial observation here is that the distance function $d_h(t,T)$ is a positive function which vanishes in the boundary and the hessian of $d_h(t,T)$ with respect to the metric $h$ is positive and bounded.\\

Using Theorem 3.18, we have
\[
n + 1 + \tilde \triangle \tilde \psi \leq  C 
\]
where $C$ is a constant independent of $T.\;$ Following estimates of Lemma 3.17,
we obtain a uniform growth control on geodesic ray.  Theorem 3.8 then follows.
%\end{proof}

\section{On the lower bound of the Calabi energy}
In this section, we will give an lower bound estimate for the
modified Calabi energy in absence of a cscK
metrc or an extremal K\"ahler metric.  Note that for algebraic manifold, the corresponding theorem is given
in [17].
\subsection{The  classic theory of Futaki-Mabuchi and A. Hwang}

Let $K = K(J)$ be a maximal compact subgroup of the automorphism group of the K\"ahler manifold and
let $\cK = \cK(J)$ be its 
Lie algebra of  gradient holomorphic vector fields in $M.\;$  According to E. Calabi,
if there is a cscK metric of CextrK metric, the cscK metric or CextrK metric must be symmetric with respect to
one of these maximal compact subgroup (up to holomorphic conjugation).   Therefore, it makes perfect sense to consider a restricted
class $\cH_K \subset \cH$ where all K\"ahler metrics are invariant under $K.\;$ For simplicity,
suppose $\omega $ is invariant with respect to action of $K.\;$ Recalled that the Lichenowicz operator
is defined as:
\[
\cL_g(f)  = f_{,\alpha \beta} d\,z^\alpha \otimes d\,z^{\beta}.
\]
where the right hand side is (2,0) component of the Hessian form of $f$ with respect to the K\"ahler metric
$g.\;$
For any metric $g\in \cH_K$,  define $Ker(\cL_g)$ as the real part of the Kernel subspace\footnote{Usually, the Kernal space can not be split as real part and imaginary part.  However, in the
case when the metric is invariant under $K(J)$,  its Lie algebra $\cK(J)$  always
corresponds to this real part of the Kernel space.} of the operator $\cL_g$ in
$C^\infty(M).\;$  It is easy to see the correspondence between $Ker(\cL_g)$ and $\cK$ in the following
formula
\[
   X = g^{\alpha \b \beta} {{\p }\over {\p w^{ \alpha}}} {{\p \theta_X}\over {\p w^{\b \beta}}} 
\] 
where
\[
X \in \cK, \theta_X \in Ker(\cL_g) \qquad  {\rm and}\;\; \displaystyle \int_M\; \theta_X \;\omega_g^n = 0.
\]
It is easy to see that such a correspondence is 1-1 as long as $g \in \cK.\;$  Futaki-Mabuchi define a
bilinear form in $\cK$ by
\[
(X,Y) = \displaystyle \int_M\; \theta_X \theta_Y \omega_g^n.
\]
Here $\theta_Y$ is the holomorphic potential of $Y.\;$
Futaki-Mabuchi proved that such a bilinear form is positive definite and well defined for $\cH_K.\;$
From the definition of $\theta_X$, it is easy to see that if $g \in \cH_K$, then $\theta_X$ is real
since
\[
L_{Im(X)} \omega_g = 0, \qquad \forall \;X\in \cK,  {\rm and}\;\: g \in \cH_K.
\]
Thus, the Futaki-Mabuchi bilinear form is positive definite.  To show it is well defined, we need
to show that it is invariant when metrics varies inside $\cH_K.\;$  Let $\omega_{g(t)} = \omega_g + t \sqrt{-1} \p \b \p \varphi \in \cH_K.\;$ Let
\[
 X = g(t)^{\alpha \b \beta} {{\p }\over {\p w^{ \alpha}}} {{\p \theta_X(t) }\over {\p w^{\b \beta}}} , \;\;{\rm and}\;\;
\]
\[
 Y = g(t)^{\alpha \b \beta} {{\p }\over {\p w^{ \alpha}}} {{\p \theta_Y(t)}\over {\p w^{\b \beta}}} .
\]
Then,
\[
\theta_X(t) = \theta_X + t X(\varphi), \qquad \theta_Y(t) = \theta_Y + t Y(\varphi),
\]
where
\[
L_{Im(X)} \varphi = L_{Im(Y)} \varphi =0.
\]
Set
\[
(X,Y)_t = \displaystyle \int_M\; \theta_X(t) \theta_Y(t) \omega_{g(t)}^n.
\]
It is straightforward to  compute
\[
\begin{array}{lcl} {d\, \over {d\, t}} (X,Y)_t  & = &  \displaystyle \int_M\; \left( \theta_Y(t) {d\, \over {d\, t}}  \theta_X(t)   + \theta_X(t) {d\, \over {d\, t}}  \theta_Y(t) + \theta_X(t) \theta_Y(t) \triangle_{g(t)} \varphi\right) \omega_{g(t)}^n\\
& = &  \displaystyle \int_M\; \left( \theta_Y(t) X(\varphi)  + \theta_X(t) Y(\varphi)\right. \\ & & \qquad \qquad \left. - g(t)^{\alpha \b \beta} {{\p  \theta_X(t)}\over {\p w^{\b \beta}}}  Y(\varphi)  {{\p \varphi}\over {\p w^\alpha}}   - g(t)^{\alpha \b \beta} {{\p  \theta_Y(t)}\over {\p w^{\b \beta}}}  X(\varphi)  {{\p \varphi}\over {\p w^\alpha}} \right) \omega_{g(t)}^n\\
& = &  \displaystyle \int_M\; \left( \theta_Y(t) X(\varphi)  + \theta_X(t) Y(\varphi) - \theta_Y(t) X(\varphi)  -\theta_X(t) Y(\varphi)\right) \omega_{g(t)}^n\\ & = & 0.
\end{array}
\]
Thus, the Futaki-Mabuchi bilinear form is well defined.  Now,  Futaki character defines a linear map
from $\cK$ to $R, $ by Rezzi representation formula, there is a unique vector field $\cX_c \in \cK$ such that
\[
  \cF_X([\omega]) = (X, \cX_c), \qquad \forall X \in \cK.
\]
Since both $\cF_X$ and the Futaki-Mabuchi form is independent of the metric, so is $\cX_c$ 
{\it a priori}.   When there is an extremal K\"ahler metric, then $\cX_c$ coincides with the complex gradient vector field of the scalar curvature function.

\begin{theo} (Hwang) The following inequality holds
\[
\displaystyle \inf_{g \in \cH_K}\; Ca(\omega_g) \geq \cF_{\cX_c}([\omega]),
\]
where the equality holds if there is an extremal K\"ahler metric in $[\omega]$.
\end{theo}
\begin{proof} Suppose $g \in \cK$ and using the $L^2$ norm with respect to $\omega_g^n$ to
decompose $ R(g) - \b R$ as 
\[
 R(g) - \b R = -  \rho  - \rho^{\perp} = -\triangle_g F, \qquad {\rm where}\;\; \rho \in Ker(\cL_g),
\]
it is easy to see that
\[
 \cX_c = g^{\alpha \b \beta} {{\p }\over {\p w^{ \alpha}}} {{\p \rho}\over {\p w^{\b \beta}}} .
\]
Thus,
\[
\begin{array}{lcl} Ca(\omega_g) & = & \displaystyle \int_M\; (R(g) - \b R)^2 \omega_g^n \\
& = &  \displaystyle \int_M\; \rho^2 \omega_g^n  +  \displaystyle \int_M\; (\rho^\perp)^2 \omega_g^n 
\\ & \geq & - \displaystyle \int_M\; \rho (R(g) - \b R) \omega_g^n  =  \displaystyle \int_M\; \rho \triangle_g F\omega_g^n  =  - \displaystyle \int_M\; \nabla \rho \cdot \nabla F \omega_g^n  \\
& = &   \displaystyle \int_M\; \cX_c (F) \omega_g^n = \cF_{\cX_c}([\omega]).  
\end{array}
\]
\end{proof}
At the time, Andrew Hwang thought the same proof could be extended to cover the non-invariant case.
Unfortunately, the Futaki-Mabuchi form is no longer positive definite and the whole argument collapsed.   Much efforts have been made by other mathematicians to bridge the gap, none are succsful.  Nonetheless, 
it is 
very interesting and also important to  
generalize the above theorem to a more general settings.

\subsection{The first derivatives of the K energy}

It is well known that the first derivatives of the K energy functional is monotonically increasing along any
smooth geodesic segment or ray.  Using the latest result of Chen-Tian \cite{chentian005},
we can show that the difference of the first derivatives of the K energy at the two ends of any
$C^{1,1}$ geodesic segment always
have a preferred sign.  This property turns out to be sufficient for our purpose.\\

For any two K\"ahler potentials $\phi_0, \phi_1\in \cal H,$ we
want to use the almost smooth solution to approximate the $C^{1,1}$
geodesic between $\phi_0$ and $\phi_1.\;$ This is an approach
first taken in \cite{chentian005}.  For any integer $l, $ consider
Drichelet problem for the HCMA equation \ref{eq:hcma0} on the
rectangle domain $\Sigma_l = [-l,l] \times [0,1]$ with boundary
value
\begin{equation}
\phi(s,0) = \phi_0, \phi(s,1) = \phi_1;\qquad \phi(\pm l,t) =
(1-t)\phi_0 + (1-t) \phi_1, \qquad (s,t) \in \Sigma_l.
\label{appl:boundarymap}
\end{equation}
We
may modify this boundary map in the four corners so that the
domain is smooth without corner. Denote the almost smooth solution
by $\phi^{(l)}: \Sigma_l \rightarrow \cal H\;$  which corresponds
to this boundary map\footnote{We may need to alter the boundary
value slightly. }.  According to \cite{chen991}, $\phi^{(l)}$ has
a uniform $C^{1,1}$ upper bound which is independent of $l.\;$ Set
\begin{equation}
  \bE^{(l)}(s,t) = \bE(\phi^{(l)}(s,t)), \qquad\forall\; (s,t) \in
  \Sigma_l.
  \label{appl:kenergyondisk}
\end{equation}
and
\begin{equation}
\bE(s,t) = \bE(\varphi(t)), \qquad\forall\; (s,t) \in S^1\times
[0,1]\times M. \label{appl:kenergyondisk1}
\end{equation}
Before we prove the main theorem, we need a convergence lemma
\begin{lem} \label{lem:approxgeode1} For any $m > 0$ fixed,  we have $\{\phi^{(l)}(s,t), l\in \bN \}$
converges uniformly to $\varphi$ in the weak $C^{1,1}$ topology.  In
particular,
\[\displaystyle \lim_{l\rightarrow \infty}\; {{\p \phi^{(l)}}\over
{\p s}} =0, \qquad \;\;\forall\;\;(s,t) \in \Sigma^{(m)}
\]
with respect to any $C^{\alpha} (0 < \alpha < 1)$ norm in
$\Sigma^{(m)}\times M.\;$
\end{lem}
\begin{proof} Note first that $\phi^{l}$ has uniform $C^{1,1}$ bound
on $\Sigma^{(m)}\times M.\;$ Thus, passing to a subsequence if
necessary, we have $\phi^{(l)}\rightarrow \b \varphi$ strongly in
$C^{1,\alpha}(\alpha\in (0,1)) $ or $W^{2,p} (p \in (0,\infty))$
in $\Sigma^{(m)}\times M$ for any $m$ fixed.   On the other hand,
$ {{\p \phi^{(l)}}\over {\p s}}$ satisfies the following equation:
\[
  \triangle_z \; {{\p \phi^{(l)}}\over {\p s}} = 0
\]
where $\triangle_z$ represents the Lapalacian operator along each
holomorphic leaf.  In other words,   $ {{\p \phi^{(l)}}\over {\p
s}}$
is a harmonic function which vanishes in main part of the boundary:
\[
  {{\p \phi^{(l)}}\over {\p s}}(s,1) = {{\p \phi^{(l)}}\over {\p s}}(s,0) =0, \qquad \forall\; s\in [-l+1,l-1].
\]
We claim that ${{\p \b \varphi}\over{\p s}} = 0.\;$ Note
that ${{\p \phi^{(l)}}\over {\p s}}\rightarrow {{\p \b
\varphi}\over{\p s}}$ in any $C^{\alpha}(\Sigma^{(m)}\times M) $
norm.  Picking any point $(z_0,x_0)\in \Sigma^{(m)}\times M, $ we
have
\[
\lim_{l\rightarrow \infty}\; {{\p \phi^{(l)}}\over {\p s}} ={{\p
\b \varphi^{(l)}}\over {\p s}}(z_0,x_0).
\]
Using this fact, we can consider the pull back function of ${{\p
\phi^{(l)}}\over {\p s}}$ on each long strip. Suppose $\cS_l$ is a
holomorphic leaf in $\Sigma_l\times M$ which passes through the
point $(z_0,x_0).\;$ Let $\jmath_l: \Sigma^{(l)} \rightarrow
\Sigma^{(l)}\times M$ be the holomorphic map associated with
this leaf $\cS_l.\; $ Let $h_l$ be the pulled back function of
${{\p \phi^{(l)}}\over {\p s}}$ under this sequence of holomorphic
maps. For simplicity, assume the preimage of $(z_0,x_0)$ is always
$(0,
{1\over 2}).\; $ \\

Now $\{h_l,l \in \bN\}$ is a sequence of bounded harmonic functions
defined in  $\Sigma^{(l)}\;$ such that it vanishes completely in
$\{0\}\times [0,1]$ and $\{1\}\times [0,1].\;$ Note that \[
\Sigma^{(1)}\subset \Sigma^{(2)}\subset \cdots \Sigma^{(l)}
\subset \cdots \subset (-\infty,\infty) \times [0,1].
\]
Since this sequence of functions vanishes in any fixed compact
subset of $\p \Sigma^{(l)}$ and the limiting domain is an infinitely
long strip, a careful argument using the  maximum principle will imply
that $h_l$ converges strongly to $0$ in any compact subset of the
infinite strip. In particular, we have
\[
  \displaystyle \lim_{l\rightarrow \infty} h_l(z_0,x_0) =0
\]
or
\[
{{\p \b \varphi}\over {\p s}}(z_0,x_0) =
\displaystyle\lim_{l\rightarrow \infty} {{\p \phi^{(l)}}\over {\p
s}}(z_0,x_0) = \displaystyle\lim_{l\rightarrow \infty}
h_l(0,{1\over 2})=0,
\]
since $(z_0,x_0)$ is an arbitrary point in  $\Sigma^{(m)}\times M$
for any $m$ fixed.  We then prove that ${{\p \phi^{(l)}}\over {\p
s}}$ converges to $0$ in any compact subdomain of
$(-\infty,\infty)\times [0,1]\times M.\;$ Since $\phi^{(l)}$
satisfies uniform $C^{1,1}$ bound,  then $\phi^{(l)}$ converges
uniformly to $\b \varphi$ in the $C^{1,\alpha}$ norm.  Using
uniqueness for the complex Monge Ampere equation,  we can show that
$\b \varphi$ is the unique $C^{1,1}$ geodesic obtained in
\cite{chen991}.
\end{proof}

  Now we are ready to prove\footnote{This lemma should be considered as a natural extension of
  chen-tian \cite{chentian005}. }
\begin{lem}\label{appl:kenergyboundbelow1} For any two K\"ahler metrics $\varphi_0, \varphi_1 \in \cH$ with $\varphi(t,\cdot)$ being the unique $C^{1,1}$ geodesic connecting these
two metrics such that $\varphi(0, x) =\varphi_0 $ and
$\varphi(1,x) = \varphi_1$,  we have
\[
  \left(d \bE\mid_{\varphi_0}, {{\p \varphi}\over {\p t}}\mid_{t=0}\right) \leq \left(d \bE\mid_{\varphi_1}, {{\p \varphi}\over {\p
  t}}\mid_{t=1}\right).
\]
\end{lem}
\begin{rem} Even though the K energy is well defined along any
$C^{1,1}$ geodesic path,  its derivative is in general not well
defined.   Thus, the evaluation of the K energy form on ${{\p \varphi}\over
{\p t}}\mid_{t=1}$ is always bigger than  it evaluation on ${{\p
\varphi}\over {\p t}}\mid_{t=0}.\;$
\end{rem}
\begin{proof}
Let ${\kappa}:(-\infty,\infty) \rightarrow {\bf R}$ be a smooth
non-negative function such that ${\kappa} \equiv 1 $ on $[ -
{1\over 2},{1\over 2}]$ and vanishes outside  of $[-{3\over
4},{3\over 4}].\;$  Set
\[
{\kappa}^{(l)} (s) = {1\over v} {\kappa}({s\over l}), \qquad {\rm
where} \;v = \int_{-\infty}^\infty\; {\kappa}(s)\;d\,s.
\]
Set
\[
f^{(l)}(t) = \displaystyle \int_{-\infty}^\infty\;
{\kappa}^{(l)}(s) {{d\,\bE}\over {d\,t}}(s,t)\;d\,s.
\]
and for any integer $m < l$, we can also define
\[
f^{(ml)}(t) = \displaystyle \int_{-\infty}^\infty\;
{\kappa}^{(m)}(s) {{d\,\bE^{(l)}}\over {d\,t}}(s,t)\;d\,s.
\]
Then
\[\begin{array}{lcl}
f^{(ml)}(0) & = & \displaystyle \int_{-\infty}^\infty\;
{\kappa}^{(m)}(s) {{d\,\bE^{(l)}}\over {d\, t}}\mid_{(s,0)}\;d\,s
\\
f^{(ml)}(1) & = &\displaystyle \int_{-\infty}^\infty\;
{\kappa}^{(m)}(s) {{d\,\bE^{(l)}}\over {d\,t}}\mid_{(s,1)} \;d\,s
%\\ {{d f^{(ml)}}\over {d\,t}}\mid_{t=0} & = & 0.
\end{array}
\]
Now
\[
\begin{array}{lcl} f^{(ml)} (1) - f^{(ml)} (0)
& = & \displaystyle \int_0^1\;  {{d f^{(ml)}}\over
{d\,t}}\;d\,t\; \\
& = & \displaystyle \int_0^1\;\displaystyle \int_{-\infty}^\infty
{\kappa}^{(m)}(s)\; {{\p^2 \bE^{(l)}}\over
{\p\,t^2}}\;d\,s\;d\,t\;
\\ & = & \displaystyle \int_0^1\;
\displaystyle \int_{-\infty}^\infty {\kappa}^{(m)}(s)\;
\triangle_{s,t} \bE^{(l)}\;d\,s\;d\,t - \displaystyle \int_0^1\;
\displaystyle \int_{-\infty}^\infty {\kappa}^{(m)}(s)\; {{\p^2
\bE^{(l)}}\over {\p\,s^2}}\;d\,s\;d\,t\\ & \geq & - \displaystyle
\int_0^1\; \displaystyle \int_{-\infty}^\infty {\kappa}^{(m)}(s)\;
{{\p^2 \bE^{(l)}}\over {\p\,s^2}}\;d\,s\;d\,t\\ & = & -
\displaystyle \int_0^1\;\displaystyle \int_{-\infty}^\infty {{d^2
{\kappa}^{(m)}(s)}\over {d\,s^2}}\; \bE^{(l)}(s,t)\;d\,s\;d\,t\;
\\ & = & -{1\over m^2} {1\over v} \displaystyle \int_0^1\;
\displaystyle \int_{-\infty}^\infty {{d^2 {\kappa}^{(m)}}\over
{d\,s^2}}\mid_{s\over m}\; \bE^{(l)}(s,t)\;d\,s\;d\,t\;
\end{array}
\]
Note that $|\bE^{(l)}(s,t)| $ has a unform bound $C.\;$ Thus,
\[
\begin{array}{lcl} {1\over m^2} {1\over
v} \mid \displaystyle \int_{-\infty}^\infty {{d^2
{\kappa}^{(m)}}\over {d\,s^2}}\mid_{s\over m}\;
\bE^{(l)}(s,t)\;d\,s \mid & \leq & {1\over m^2} {1\over v}
\displaystyle \int_{-\infty}^\infty |{{d^2 {\kappa}^{(m)}}\over
{d\,s^2}}\mid_{s\over m}|\; \bE^{(l)}(s,t)\;d\,s
\\& \leq & {C\over {v\; m^2}}\;\displaystyle \int_{-\infty}^\infty |{{d^2
{\kappa}^{(m)}}\over {d\,s^2}}\mid_{s\over m}|\; \;d\,s \\
& = & {C\over {v\; m}}\;\displaystyle \int_{-\infty}^\infty |{{d^2
{\kappa}^{(m)}}\over {d\,s^2}}\mid_{s}|\; \;d\,s = {C\over {v\;
m}}\;\displaystyle \int_{-1}^1  |{{d^2 {\kappa}^{(m)}}\over
{d\,s^2}}\mid_{s}|\; \;d\,s\\
& \leq & {C\over m}
\end{array}
\]
for some uniform constant $C.\;$ Therefore, we have
\begin{eqnarray}
f^{(ml)} (1) - f^{(ml)} (0) & \geq & -{1\over m^2} {1\over v}
\displaystyle \int_0^1\; \displaystyle \int_{-\infty}^\infty {{d^2
\kappa}\over {d\,s^2}}\mid_{s\over m}\;
\bE^{(l)}(s,t)\;d\,s\;d\,t\;\nonumber
\\
&\geq & - \displaystyle \int_0^1\;{1\over m^2} {1\over v}\;\mid
\displaystyle \int_{-\infty}^\infty {{d^2 {\kappa}}\over
{d\,s^2}}\mid_{s\over m}\; \bE^{(l)}(s,t)\;d\,s\mid
\;d\,t\nonumber
\\&\geq &- \displaystyle \int_0^1\;  {C\over
m} \;d\,t = -{C\over {2 m}}. \label{eq:firstderivativeofKenergy}
\end{eqnarray}
Recall that
\[\begin{array}{lcl}
f^{(ml)}(1) & = & \displaystyle \int_{-\infty}^\infty\;
{\kappa}^{(m)}(s) {{d\,\bE^{(l)}}\over
{d\,t}}(s,1)\;d\,s\\
& =&  \displaystyle \int_{-\infty}^\infty\; {\kappa}^{(m)}(s)
\int_M\; (\underline{R}-R(\varphi_1)) {{\p\phi^{(l)}}\over {\p t}}
\omega_{\varphi_1}^n\;d\,s .\end{array}
\]

For any fixed $m$, by Lemma \ref{lem:approxgeode1},
$\phi^{(l)}$ uniformly converges to $\varphi$ in
$\Sigma^{(m)}\times M.\;$ In particular, ${{\p \phi^{(l)}}\over
{\p t}}$ converges strongly in $C^{\alpha}$ norm to ${{\p
\varphi}\over {\p t}}$ in $\Sigma^{(m)}\times M.\;$ Thus, fixing  $m$
 and letting $l \rightarrow \infty,$ we have
\[
\begin{array} {lcl} \displaystyle \lim_{l \rightarrow \infty} f^{(ml)}(1) &
= & \displaystyle \int_{-\infty}^\infty\; {\kappa}^{(m)}(s)
\displaystyle \int_M\; (\underline{R}-R(\varphi_1))
{{\p\varphi}\over {\p t}}\mid_{t=1}\;\omega_{\varphi_1}^n\;d\,s\\
& = & \displaystyle \lim_{l \rightarrow \infty} \displaystyle
\int_{-\infty}^\infty\; {\kappa}^{(m)}(s) \displaystyle \int_M\;
(\underline{R}-R(\varphi_1))
{{\p\phi^{(l)}}\over {\p t}}\mid_{t=1} \omega_{\varphi_1}^n \;d\,s \\
& = & \displaystyle \int_{-\infty}^\infty\; {\kappa}^{(m)}(s)
\displaystyle \int_M\; (\underline{R}-R(\varphi_1))
{{\p\varphi}\over {\p t}}\mid_{t=1} \omega_{\varphi_1}^n \;d\,s\\
& =  & \displaystyle \int_M\; (\underline{R}-R(\varphi_1))
{{\p\varphi}\over {\p t}}\mid_{t=1} \omega_{\varphi_1}^n.
\end{array}\]

Similarly, we can prove
\[
\displaystyle \lim_{l \rightarrow \infty} f^{(ml)}(0) =
\displaystyle \int_M\; (\underline{R}-R(\varphi_1))
{{\p\varphi}\over {\p t}}\mid_{t=0} \omega_{\varphi_0}^n.
\]
Plugging this into inequality \ref{eq:firstderivativeofKenergy},
we have
\[
(d\, \bE, {{\p \varphi}\over {\p t}}\mid_{t=1})_{\varphi_1} - (d\,
\bE, {{\p \varphi}\over {\p t}}\mid_{t=0})_{\varphi_0} \geq -
{C\over m}.
\]
As $m \rightarrow \infty$, we have
\[
(d\, \bE, {{\p \varphi}\over {\p t}}\mid_{t=1})_{\varphi_1} \geq
(d\, \bE, {{\p \varphi}\over {\p t}}\mid_{t=0})_{\varphi_0}.
\]
The Lemma is
then proved.
\end{proof}
\subsection{The lower bound of the Calabi energy}
Note that the first derivative of the K energy
functional is always non-decreasing along a geodesic ray.  Thus,
the $\yen$ invariant is always well defined along any relative $C^{1,1}$ geodesic ray.  Now
we are ready to prove Theorem \ref{th:calabienergybound0}.
\begin{proof}  Suppose $\rho:[0,\infty)\rightarrow\infty$ is a
geodesic ray parametrized by arc length such that
\[
  \displaystyle \lim_{t\rightarrow \infty}(d\,\bE, {{\p\rho}\over
  {\p t}})_{\rho(t)} < -\infty.
\]
For any K\"ahler potential $\varphi_0\in \cH$, consider the unique
$C^{1,1}$ geodesic connecting $\varphi_0$ to $\rho(l).\;$  Represent 
this geodesic by $\psi_l:[0, L_l] \rightarrow \cH\;$ and parametrize by arc length.  Applying the
preceding Theorem,  we have
\[\begin{array}{lcl} &&
(d\,\bE, {{\p\psi_l}\over {\p s}}\mid_{s=0})_{\varphi_0} \\& \leq & (d\,\bE, {{\p\psi_l}\over {\p
s}}\mid_{s=L_l})_{\rho(l)} \\ & = & (d\,\bE, {{\p\psi_l}\over {\p
s}}\mid_{s=L_l} - {{\p \rho}\over {\p t}}\mid_{t=l} )_{\rho(l)} + (d\,\bE,{{\p \rho}\over {\p t}}\mid_{t=l}
)_{\rho(l)}\\ & \leq & \left(\int_{M}\; (R(\rho(l)) - \underline{R})^2
\omega_{\rho}^n\right)^{1\over 2} \cdot \left(  \int_M\;( {{\p\rho}\over {\p
t}}\mid_{s=L_l} -  {{\p\psi_l}\over {\p
t}}\mid_{s=L_l} )^2 \omega_{\rho(l)}^n\right)^{1\over 2}  + (d\,\bE,{{\p \rho}\over {\p t}}\mid_{t=l})_{\rho(l)}\\
& \leq  & C ( 1 - ( {{\p\rho}\over {\p
t}}\mid_{t=l},\; {{\p \psi_l}\over {\p s}}_{s=L_l})_{\rho(l)})^{1\over 2}+ (d\,\bE,{{\p \rho}\over {\p
t}}\mid_{t=l})_{\rho(l)}
\end{array}
\]
Since $\cH$ is a non-positively curved manifold in the sense of
Alexander,  we have
\[
( {{\p\rho}\over {\p
t}}\mid_{l},\; {{\p \psi_l}\over {\p s}}_{s=L_l} )_{\rho(l)}  \rightarrow
1
\]
as $l \rightarrow \infty.\;$  Then,  we have
\[
\begin{array}{lcl} \yen(\rho) &= &\displaystyle \liminf_{l\rightarrow \infty}\; (d\,\bE,{{\p
\rho}\over {\p
t}}\mid_{l})_{\rho(l)} 
\leq 
(d\,\bE, {{\p\psi_l}\over {\p s}}\mid_{s=0})_{\varphi_0} \\
&\leq &  \left(\displaystyle \int_M (R(\omega_{\varphi_0}) - \b R)^2 \omega_{\varphi_0}^n\right)^{1\over 2}
\left( \displaystyle \int_M\; ({{\p\psi_l}\over {\p s}}\mid_{s=0})^2\; \omega_{\varphi_0}^n\right)^{1\over 2}\\
& = & (Ca(\omega_{\varphi_0}))^{1\over 2}.
\end{array}
\]
In other words, we have
\[
Ca(\omega_{\varphi_0}) \geq \yen(\rho)^2.
\]
Our theorem follows from here directly!
\end{proof}
Now we are ready to prove Theorem \ref{th:calabienergybound3}.
\begin{proof} Let $\cX_c$ to be the {\it a priori} extremal vector field. Suppose $g \in \cH_K.\;$
Suppose that $\omega_{\rho(t)} (t\in (-\infty, \infty))$ is the one paramter family of K\"ahler metrics
generated by pulling the K\"ahler metrics $\omega_g$ in the direction of $Re(X).\;$  It is straightforward
to check that $\rho(t)$ satisfies the geodesic equation and
\[
{{d\, E(\omega_{\rho(t)})}\over {d\, t}} =\pm \cF_{\cX_c}([\omega]).
\]
Select one direction so that
\[
{{d\, E(\omega_{\rho(t)})}\over {d\, t}} = - \cF_{\cX_c} ([\omega]).
\]
Note that the length element of this geodesic line
is
\[\begin{array}{lcl}
\int_M\; \left({{\p \rho}\over {\p t}}\right)^2 \omega_{\rho(t)}^n & = & (\theta_{\cX_c}, \theta_{\cX_c})
  =  - (\theta_{cX_c}, R(\rho) - \b R)\\ &  = & - \int_M\; {{\p \rho}\over {\p t}} (R(\rho(t)) - \b R) \omega_{\rho(t)}^n  = \cF_{\cX_c}.
\end{array}
\]
Now, if we re-paramterize by arc length, then the $\yen$ invariant along this
geodesic line satisies
\[
\yen(\rho)^2  = \cF_{\cX_c}([\omega]).
\]
Our theorem then follows from Theorem \ref{th:calabienergybound0}.  

\end{proof}

\section{On the lower bound of the geodesic distance}
Let us   prove Theorem 1.2 first.
\begin{proof}  We follow the notations of Subsection 5.2.    Set
\[
  E^{(ml)}(t)  =  \displaystyle \int_{-\infty}^\infty\; k^{(m)}(s) \; E^{(l)}(s,t) \;d\,s, \qquad \forall m \leq \l \in \NN.
\]
Then,
\[
  E^{(ml)}(0) = E(\varphi_0),  \; E^{(ml)}(1) = E(\varphi_1)
\]
and
\[
 {{d\, E^{(ml)}}\over {d\,t}}(t) = f^{(ml)}(t), \qquad \forall\; t\in [0,1].
\]
Following the same calculation as in subsection 5.2, for any $ 0\leq t_1 < t_2 \leq 1$
\begin{eqnarray}
f^{(ml)} (t_2) - f^{(ml)} (t_1) & \geq & -{1\over m^2} {1\over v}
\displaystyle \int_{t_1}^{t_2}\; \displaystyle \int_{-\infty}^\infty {{d^2
\kappa}\over {d\,s^2}}\mid_{s\over m}\;
\bE^{(l)}(s,t)\;d\,s\;d\,t\;\nonumber
\\
&\geq & - \displaystyle \int_{t_1}^{t_2}\;{1\over m^2} {1\over v}\;\mid
\displaystyle \int_{-\infty}^\infty {{d^2 {\kappa}}\over
{d\,s^2}}\mid_{s\over m}\; \bE^{(l)}(s,t)\;d\,s\mid
\;d\,t\nonumber
\\&\geq &- \displaystyle \int_{t_1}^{t_2}\;  {C\over
m} \;d\,t = -{C\over {2 m}}. \label{eq:firstderivativeofKenergy1}
\end{eqnarray}
%\end{proof}

Thus,
\[
  f^{(ml)}(0) - {C\over {2 m}} \leq  f^{(ml)}(t) \leq f^{(ml)}(1) + {C\over {2 m}}
\]
Therefore, 
\[
\begin{array}{lcl}  E(\varphi_1) - E(\varphi_0) & = &  E^{(ml)}(1) - E^{(ml)}(0)\\
& = & \displaystyle \int_0^1\;  {{d\, E^{(ml)}}\over {d\,t}}(t)  \; d\,t =  \displaystyle \int_0^1\; f^{(ml)}(t) \; d\,t\\
& \leq  &  \displaystyle \int_0^1\;\left( f^{(ml)}(1) + {C\over {2m}}\right)\; d\,t\\
& = &  \displaystyle \int_{-\infty}^\infty\; {\kappa}^{(m)}(s)
\int_M\; (\underline{R}-R(\varphi_1)) {{\p\phi^{(l)}}\over {\p t}}
\omega_{\varphi_1}^n\;d\,s  + {C\over {2m}}
\end{array}
\]
As before, let $l \rightarrow \infty$, so $\phi^{(l)}(s,t) $ converges to the geodesic $\varphi(t)$
strongly in $C^{1,\alpha}$ norm.   Then, letting $m \rightarrow \infty$, we have
\[\begin{array}{lcl}
E(\varphi_1) - E(\varphi_0) & \leq &   \displaystyle 
\int_M\; (\underline{R}-R(\varphi_1)) {{\p\varphi}\over {\p t}}
\omega_{\varphi_1}^n\\
& \leq & \sqrt{Ca(\varphi_1)} \sqrt{\displaystyle \int_M\;  \left( {{\p\varphi}\over {\p t}}\right)^2 \omega_{\varphi_1}^n} \\
&= & \sqrt{Ca(\varphi_1)} \cdot d(\varphi_0, \varphi_1).
\end{array}
\]
\end{proof}
Corollary 1.3 follows from this theorem since the $|\varphi|_{\infty}$ bound will imply the geodesic distance of $\varphi$ to a fixed K\"ahler potential is bounded.\\

Before we prove Theorem 1.4, we need to prove a proposition first.
\begin{prop}\cite{chentian0}  Let $Ric(\omega_\varphi) \geq - C_1$ then there is a uniform constant $C$ such
that :\beg
\inf_M \log {\omega_\varphi^n \over \omega^n}  \geq-4C\exp\left(2+2\int_M\log\frac{\omega^n_{\varphi}}{\omega^n}\omega^n_{\varphi}\right).
\ee
If $C_1 = 0, $ then $C$ is dimensional constant. Otherwise, $C$ depends on $C_1$ and $|\varphi|_{L^\infty}\;$ or $\sup \varphi_-  + \int_M\; \varphi_+ \omega^n.\;$
\end{prop}
\begin{proof} Set \[
F =  \log {\omega_\varphi^n \over \omega^n} .
\]  The lower Ricci curvature bound implies  that
\[
  Ric(\omega) -  i \p \b \p F \geq - C_1  \omega_\varphi. 
\]
Taking the trace of both sides, we have
\[
\triangle ( F  - C_1 \varphi) \leq C_2
\]
for some constant $C_2.\;$\\

Choose a constant $c$ such that
\beg\int_M\log\frac{\omega^n_{\varphi}}{\omega^n}\omega^n_{\varphi}\leq
c,\ee In a fixed K\"ahler class, we have
\[
\int_M\omega^n=Vol(M)=1.
\]
Put $\epsilon$ to be $\exp(-2-2c)$. Observe that
\[
\log\frac{\omega^n_{\varphi}}{\omega^n}\omega^n_{\varphi}\geq-e^{-1}\omega^n,
\]
We have \beg c&&\geq
\left(\int_{\epsilon\omega^n_{\varphi}>\omega^n}+\int_{\epsilon\omega^n_{\varphi}\leq\omega^n}\right)
\left(\log\frac{\omega^n_{\varphi}}{\omega^n}\omega^n_{\varphi}\right)\\
&&\geq\int_{\epsilon\omega^n_{\varphi}>\omega^n}\left(\log\frac{1}{\epsilon}\right)\omega^n_{\varphi}
+\int_{\epsilon\omega^n_{\varphi}\leq\omega^n}\left(-e^{-1}\omega^n\right)\\
&&>2(1+c)\int_{\epsilon\omega^n_{\varphi}>\omega^n}\omega^n_{\varphi}-1.\ee
It follows that
\beg\int_{\epsilon\omega^n_{\varphi}>\omega^n}<\frac{1}{2},\ee and
consequently, \beg
\int_{\omega^n\leq4\omega^n_{\varphi}}\omega^n\geq\epsilon\int_{\frac{\epsilon}{4}\omega^n\leq
\epsilon\omega^n_{\varphi}\leq\omega^n}\omega^n_{\varphi}\geq\frac{1}{4},\ee
because we know
\[\int_{\omega^n\leq4\omega^n_{\varphi}}\omega^n_{\varphi}>\frac{3}{4}
\]
and
\[\int_{\omega^n\leq\epsilon\omega^n_{\varphi}}>\frac{1}{2}.\]
Now by Green's formula, we have \beg ( F - C_1 \varphi)(p)=-\int_M\; G(p,q) \triangle
(F - C_1 \varphi) \omega^n(q)+\int_M\; (F - C_1\varphi) \omega^n,\ee where $G(p,q)\geq 0$ is a Green
function of $\omega$. If either $|\varphi|_{L^\infty}$ is bounded,  or
\[
\sup \varphi_- \leq C, \qquad {\rm and}\; \; \int_M\; \varphi_+\; \omega^n \leq C,
\]
then
\beg\inf_M F & \geq &  \inf_MF\int_{\omega^n\geq4\omega^n_{\varphi}}\omega^n+\int_{\omega^n<4\omega^n_{\varphi}}F\omega^n - C_1 \sup \varphi + C_1 \int_M (-\varphi_-) \omega^n   \\&\geq& \int_M F\omega^n-C\\
& \geq & \inf_MF\int_{\omega^n\geq4\omega^n_{\varphi}}\omega^n+\int_{\omega^n<4\omega^n_{\varphi}}F\omega^n-C\\
&\geq & \inf_MF\int_{\omega^n\geq4\omega^n_{\varphi}}\omega^n-\log4\int_{\omega^n<4\omega^n_{\varphi}}\omega^n-C\\
&\geq &\left(1-\frac{\epsilon}{4}\right)\inf_MF-C,\ee where we can
assume $\inf_M F<0.$ Therefore, we have
\[\inf_MF\geq-4C\exp(2+2c).\] By the way we choose the constant
$c$ in the beginning of the proof, proposition is proved.
\end{proof}

One more lemma is needed.
\begin{lem} Let $\omega'\in [\omega], $ and suppose  $(M, \omega')$ has a bounded
Soblev constant.   For any $\psi \in C^\infty(M), $ and $\; \omega' + i \p \b \p \psi > 0$,
then $\psi_+ = \max (\psi, 0)$ is uniformly bounded if $\displaystyle \; \int_M\; \psi_+^p \omega'^n < C\;$
for any $p > 1.\;$
\end{lem}
The proof is well known to the experts and we will include it here for the convenience of readers.
\begin{proof}  Without loss of generality, may assume $\psi \geq 1$ for simplicity.   We start from
\[
 n + \triangle' \psi > 0, 
\]
here $\triangle'$ is the Laplacian operator of $\omega'.\;$  For any $p \geq 1,$ we have
\[\begin{array}{lcl}
  \int_M\; n \cdot \psi^p  \omega'^n & \geq & p\int_M\; |\nabla \psi|^2 \psi^{p-1} \omega'^n\\
   & = &
   p\int_M\; |\nabla \psi^{{p+1}\over 2}|^2 \psi^{p-1} \omega'^n\\
   & = & {{4 p}\over {(p+1)^2}} \int_M |\nabla \psi^{{p+1}\over 2}|^2 \omega'^n. \end{array}
\]
Since $\omega'$ has a uniform Soblev constant, we have
\[\begin{array}{lcl} c_{Sob} \left(\int_M\; \psi^{{{p+1}\over 2} \cdot {{2m}\over {m-2}}}\right)^{{m-2}\over {m}} & \leq & \int_M\; |\nabla \psi^{{p+1}\over 2}|^2  + (\psi^{{p+1}\over 2})^2\\
& \leq & {(p+1)^2 \over { 4 p}}  n \int_M\; \psi^p \omega'^n + \int_M\; \psi^{p+1}\omega'n
\end{array}\]
%\end{proof}
Thus, there is a uniform constant $C$ which is independent of $p$ such that
\[
 \left(\int_M\; \psi^{ (p+1)\cdot {{m}\over {m-2}}}\right)^{{m-2}\over {m}}  \leq C (p+1)  \int_M\; \psi^{p+1}\omega'^n.
\]
Set 
\[
p_1  = q > 0,  p_2 = p_1 {m\over {m-2}}, \cdots p_{j+1} = p_j \cdot   {m\over {m-2}}, \cdots
\]
Then,
\[
\|\psi\|_{p_{j+1}} \leq C^{1\over p_j} {p_j}^{1\over p_j} \|\psi\|_{p_j}, \qquad, \forall j = 1,2, \cdots.
\]
In other words,
\[
\|\psi\|_{p_{j+1}} \leq \|\psi\|_{p_1} \cdot C^{\sum_{k=1}^j\left( {1\over p_j} + {1\ p_j} \log p_j\right)}\|\psi\|_{p_1}.
\]
Let $p_1 = 1,$ and $j\rightarrow \infty$, then
\[
\|\psi\|_{L^\infty} \leq C \|\psi\|_{L^2}.
\]
\end{proof}
Now we are ready to prove Theorem 1.4.
\begin{proof}  Suppose $\varphi \in \cH$ is a K\"ahler potential such that
\begin{enumerate}
\item $Ric(\omega_\varphi) $ is bounded from below;
\item The diameter of $(M, \omega_\varphi)$ is bounded from above;
\item the geodesic distance $d(0, \varphi)$ is bounded from above.
\end{enumerate}
The first two conditions implies that there is a uniform Soblev and Poincare constants for $(M, \omega_\varphi).\;$  In the proof here, ``C" represents a generic constant.\\

Now normalize $\varphi $ by a constant necessary so we have
\[
I(\varphi) = 0.
\]
According to a theorem in \cite{chen991},  we have
\begin{equation}
 \displaystyle \max \left( \displaystyle \int_M \; \varphi_- \omega^n, \displaystyle \int_M \; \varphi_+ \omega_\varphi^n,\right) \leq d(0, \varphi) \leq C. \label{eq:geodesicdistance0}
\end{equation}
Here $\varphi_+, \varphi_-$ are positive and negative part of $\varphi$ respectively.\\

Set the K energy of $\omega_0$ being $0.\;$ Theorem 1.2 implies that 
\[
\bE(\varphi) \leq \bE(0)  + \sqrt{Ca(\varphi)} d(0,\varphi) \leq C.
\]
If the K energy functional is quasi-proper  ( c.f., the detailed expression of the K energy
functional \cite{chen00}),  we obtain
\[
\int_M\; \log {{\omega_\varphi^n}\over \omega^n} \omega_\varphi^n \leq C.
\]
According G. Tian, there is a positive constant $\alpha > 0$ which depends only
on the polarization $(M, [\omega])$ such that for any $\varphi \in \cH$, we have
\[
\int_M\; e^{-\alpha (\varphi -\sup \varphi)} \omega^n \leq C.
\]
Or 
\[
\int_M\; e^{-\alpha (\varphi-\sup\varphi) - \log {{\omega_\varphi^n}\over \omega^n}} \omega_\varphi^n \leq C.
\]
Consequently, we have
\[
\int_M\; -\alpha (\varphi - \sup \varphi) - \log  {{\omega_\varphi^n}\over \omega^n} \omega_\varphi^n \leq C.
\]
Therefore, 
\[\begin{array}{lcl}
\alpha \sup \varphi & \leq & \alpha \int_M\; \varphi \omega_\varphi^n  + \int_M \;  \log  {{\omega_\varphi^n}\over \omega^n} \omega_\varphi^n\\
& \leq &  \alpha \int_M\; \varphi \omega_\varphi^n  + \int_M \;  \log  {{\omega_\varphi^n}\over \omega^n} \omega_\varphi^n\\
& \leq & C.
\end{array}
\]
it follows,
\begin{equation}
\displaystyle \int_M\; |\varphi| \omega^n  \leq C. \label{eq:c0estimate}
\end{equation}
%Using the fact that Ricci is bounded below,  by appealing to Proposition 5.1, we have 
%\[
%\log {\omega_\varphi^n \over \omega^n} \geq -C.
%\]  
By the detailed expreesion of $I(\varphi)$ (cf. equation \ref{eq:Ifunctional1}), we have
\begin{equation}
J(\varphi) \leq C.\label{eq:jbound}
\end{equation}

Alternatively, when the K energy is proper, we can obtain estimate \ref{eq:c0estimate} 
and \ref{eq:jbound}
as well.  Note that if the K energy functional $\bE$ is proper in $\cH$, we have
\begin{equation}
0 \leq J(\varphi) \leq C. \label{eq:jbound1}
\end{equation}
Again, from the detailed expression of $I(\varphi)$  (cf. equation \ref{eq:Ifunctional1}), we have
\[
|\displaystyle \int_M\; \varphi \omega^n |  \leq C.
\]
Comparing to estimate \ref{eq:geodesicdistance0}, we have
\begin{equation}
 \int_M \varphi_+ \omega^n \leq C \label{eq:l1bound1}
\end{equation}
and
\begin{equation}
\int_M\; \varphi_- \omega_\varphi^n \leq C. \label{eq:l1bound2}
\end{equation}
or
\begin{equation}
\displaystyle \int_M\; |\varphi| \omega^n  \leq C. \label{eq:c0estimate1}
\end{equation}

Since $Ric(\omega_\varphi) \geq -C$ and the diameter is bounded from below, we have a uniform 
Poincare constant for $(M, \omega_\varphi).\;$
Using the Poincare inequality, we have
\[
\int_M \varphi^2 \omega^n  + \int_M \; \varphi^2 \omega_\varphi^n \leq C\left( J(\varphi) + \left( \int_M\; |\varphi| \omega_\varphi^n \right)^2\right) \leq C.
\]

Recall that 
\[
  n + \triangle \varphi \geq 0.
\]
Using Moser iteration, and the $J$ functional bound \ref{eq:jbound},
we obtain
\[
0 \leq \varphi_+  \leq C.
\]
Recall that
\[
 n + \triangle_\varphi (-\varphi) \geq 0.
\]
By the assumption that the Soblev constant of $(M, \omega_\varphi)$ is bounded and the $L^2$ norm is bounded above, appealing to Lemma 5.2, gives us 
\[
0  \leq \varphi_- \leq C.
\]
In other words, we have
\[
 |\varphi|_{L^\infty} \leq C.
\]

To derive an upper bound on the volume form,  first  note that $Ric(\omega_\varphi)$ is bounded from above, thus
\begin{equation}
 \triangle \log {\omega_\varphi^n \over \omega^n} + C_2 \varphi  \geq -C
 \label{eq:ricciupperbound}
\end{equation}
for some constant $C_2, C.\;$   Thus
\[
\begin{array}{lcl} && \left( \log {\omega_\varphi^n \over \omega^n} + C_2 \varphi\right)(x) \\& = & - \int_M\; G(x,y) \left( \log {\omega_\varphi^n \over \omega^n} + C_2 \varphi\right) \omega^n + \int_M\; \left( \log {\omega_\varphi^n \over \omega^n} + C_2 \varphi\right) \omega^n
\\ & \leq & + C + \log \int_M \; {\omega_\varphi^n\over \omega^n} \omega^n + \int_M \varphi \omega^n
\\ &  \leq & C + \int_M\; \varphi \omega^n. \end{array}
\]
Using the fact that $|\varphi|_{L^\infty} $ is bounded, we have
\[
\log {\omega_\varphi^n \over \omega^n} \leq C
\]
for some constant $C.\;$\\

To prove the metric is equivalent, we follow Yau's proof of the Calabi conjecture (cf. 
\cite{chenhe05}).  Folllowing notations in  Subsection 3.2.2,  set \[
u = \exp(-\lambda \varphi)(n+\triangle
\varphi), \qquad {\rm and}\;\;  F = \log {\omega_\varphi^n \over \omega^n}.
\]
At the maximal point $p$ of the function $u$, similar to inequality \ref{eq:estimate1}, we have
 \beg
\triangle_{\varphi}\left\{\exp(-\lambda\varphi)(n+\triangle
\varphi)\right\}(p)\leq0.\ee At the point $p$, we have \beg0&\geq&
\triangle F-n^2\inf_{i\neq
l}R_{i\bar{i}l\bar{l}}- \lambda n(n+\triangle \varphi)\\
&&\quad+\left(\lambda+\inf_{i\neq
l}R_{i\bar{i}l\bar{l}}\right)\exp\left\{\frac{-F}{n-1}\right\}(n+\triangle
\varphi)^{\frac{n}{n-1}} .\ee

Applying inequality \ref{eq:ricciupperbound},  we have  \beg0&\geq& -C\triangle
\varphi-C_2-n^2\inf_{i\neq
l}R_{i\bar{i}l\bar{l}}-\lambda n(n+\triangle \varphi)\\
&&\quad+\left(\lambda +\inf_{i\neq
l}R_{i\bar{i}l\bar{l}}\right)\exp\left\{\frac{-F}{n-1}\right\}(n+\triangle
\varphi)^{\frac{n}{n-1}} . \ee 
Since we already have control of  $F$ from both above and below here,  we can choose  $\lambda$ large enough, to imply that $ (n+\triangle\varphi)(p)$
is uniformly bounded from above.   Therefore, 
\[
u = \exp(- \lambda \varphi)(n+\triangle
\varphi), \qquad {\rm and}\;\;  F = \log {\omega_\varphi^n \over \omega^n}
\] is uniformly bounded from above (since $|\varphi|_{L^\infty}|$ is uniformly bounded).
Thus 
 \[ 0<n+\triangle
\varphi \leq  C\] Thus,  $\omega_\varphi$ is uniformly
equivalent to $\omega.\;$  Consequently,  $|\triangle F|$ is uniformly 
bounded.  Thus, the metric is uniformly $C^{1,\alpha}$ bounded for any $\alpha \in (0,1)$.

\end{proof}

%\bibliography{testa_bryant}
\end{document}